\documentclass{article}%
\usepackage{graphicx}
\usepackage{amsmath}
\usepackage{amsfonts}
\usepackage{amssymb}%
\numberwithin{equation}{section}
\newtheorem{theorem}{Theorem}

\newtheorem{corollary}[theorem]{Corollary}

\newtheorem{definition}[theorem]{Definition}

\newtheorem{lemma}[theorem]{Lemma}

\newtheorem{proposition}[theorem]{Proposition}
\newtheorem{remark}[theorem]{Remark}

\newenvironment{proof}[1][Proof]{\textbf{#1.} }{\ \rule{0.5em}{0.5em}}
\begin{document}
\bigskip

\bigskip

\begin{center}
\textbf{OPTIMAL GRADIENT ESTIMATES AND ASYMPTOTIC BEHAVIOUR FOR THE
VLASOV-POISSON SYSTEM WITH SMALL INITIAL DATA.}

\bigskip

Hyung Ju Hwang\footnote{School of Mathematics, Trinity College, Dublin 2,
Ireland.},\ Alan D. Rendall\footnote{Max Planck Institute for Gravitational
Physics, Albert Einstein Institute, Am M\"{u}hlenberg 1, 14476 Potsdam,
Germany.},\ Juan J. L. Vel\'{a}zquez\footnote{Departamento de Matem\'{a}tica
Aplicada, Universidad Complutense, Madrid 28040, Spain.}

\bigskip
\end{center}

\bigskip

\begin{abstract}

The Vlasov-Poisson system describes interacting systems of 
collisionless particles. For solutions with small initial data in
three dimensions it is known that the spatial density of particles 
decays like $t^{-3}$ at late times. In this paper this statement
is refined to show that each derivative of the density which is taken 
leads to an extra power of decay so that in $N$ dimensions for $N\ge 3$
the derivative of the density of order $k$ decays like $t^{-N-k}$.
An asymptotic formula for the solution at late times is also obtained.
\end{abstract}

\section{INTRODUCTION}

\bigskip

\bigskip

The Vlasov-Poisson system provides a statistical description of the dynamics
of a large number of particles which are acted on by a force field which they
generate collectively. One class of applications of this system is in plasma
physics where the force is electrostatic and the particles are electrons or
ions \cite{goldston}. Another is in stellar dynamics where stars play the role
of particles. The particle treatment is justified in models of galaxies where
the distance between stars is much larger than their diameters. In this case
the force is gravitational \cite{binney}. The equations in these two cases
only differ by a sign and a lot of the mathematical theory works in exactly
the same way for both. This applies in particular to the results of this
paper. For surveys of results on the Vlasov-Poisson and related systems see
\cite{glassey1} and \cite{andreasson}.

The distribution function $f$ of the particles satisfies the Vlasov equation
while the potential $\phi$ for the field satisfies the Poisson equation. The
function $f$ depends on time $t$, the spatial point $x\in\mathbb{R}^{3}$ and
the velocity $v\in\mathbb{R}^{3}$. It is natural to pose an initial value
problem with $f$ being prescribed at $t=0$. For an initial datum which is
$C^{1}$ and has compact support it is known that there exists a unique
corresponding $C^{1}$ solution, globally in time \cite{pfaffelmoser},
\cite{lions}. The support of $f$ is compact at each fixed time $t$ and an
important diagnostic quantity is $P(t)$, the supremum of $|v|$ over the
support of $f$ at time $t$. Estimates are known for $P(t)$ \cite{horst93} and
these imply estimates for $\Vert\rho(t)\Vert_{L^{\infty}}$ where $\rho$, the
spatial density of particles, is given by $\rho(t,x)=\int f(t,x,v)dv$. We have
$P(t)\leq C(1+t)\log(2+t)$. Unfortunately these
estimates seem far from optimal. They are the same for the plasma physics and
stellar dynamics cases. Intuitively it is to be expected that the optimal
estimates differ in these two cases. In the stellar dynamics case there exist
time-independent solutions so that $\Vert\rho(t)\Vert_{L^{\infty}}$ does not
decay in general. In the plasma physics case decay estimates for integral
norms of $\rho$ are established in \cite{ir} and \cite{perthame}. The pointwise
estimates can also be improved to give a bound for $P(t)$ of the form 
$C(1+t)^{2/3}$ \cite{rein98}.

It is possible to consider the analogue of the Vlasov-Poisson system in higher
dimensions. It is, however, known that global existence fails in four space
dimensions \cite{horst82}. An explicit example of singularity formation and
information on the asymptotics of solutions near a singularity were obtained
in \cite{lemou}.

There is a case where much more is known about the long-time asymptotics of
solutions of the Vlasov-Poisson system, namely that of small initial data. The
first global existence theorem for that case due to Bardos and Degond
\cite{bardos} naturally comes with decay estimates. They show that
\[
\Vert\rho(t)\Vert_{L^{\infty}}\leq C(1+t)^{-3}%
\]
If the data are sufficiently differentiable then the same techniques should
lead to estimates of the form
\[
\Vert D^{k}\rho(t)\Vert_{L^{\infty}}\leq C(1+t)^{-3}%
\]
but they do not give more. In this paper we apply new techniques to this
problem to obtain estimates of the form
\[
\Vert D^{k}\rho(t)\Vert_{L^{\infty}}\leq C(1+t)^{-3-k}%
\]
for solutions with small initial data. Furthermore, we obtain asymptotic
expansions for these solutions.

Note that there are a number of generalizations of the results of
\cite{bardos} in the literature. The fully relativistic generalization of the
plasma physics problem is given by the Vlasov-Maxwell system. An analogue of
the result of \cite{bardos} in that case was proved in \cite{glassey}. In the
stellar dynamics problem the fully relativistic generalization is the
Einstein-Vlasov system \cite{rendall04} which is much more complicated. A
small data global existence theorem in the spherically symmetric case was
obtained in \cite{rr92}. A related system which is physically incorrect but
mathematically interesting is the Vlasov-Nordstr\"{o}m system for which there
is a global existence theorem \cite{calogero}. Surprisingly it seems that no
analogue of the asymptotic result of \cite{bardos} has been proved for this
system. There are generalizations of the results for solutions of the
Vlasov-Poisson and Vlasov-Maxwell systems with small data to almost
spherically symmetric data \cite{schaeffer87},\cite{rein90}. There are also
results for solutions of the Vlasov-Poisson system with non-standard boundary
conditions which are relevant to cosmology \cite{rr94}, \cite{rein97}. Global
existence has also been proved for some cosmological solutions of the
Einstein-Vlasov system with symmetry. See for instance \cite{tchapnda}. It
would be interesting to extend the results of this paper to some of the cases
mentioned in this paragraph.

This paper was motivated by the wish to prove a small data global existence
theorem for the Einstein-Vlasov system which does not require any symmetry
assumptions. To understand the difficulty of this problem note first that even
the vacuum Einstein equations, from the present point of view the
Einstein-Vlasov system with $f=0$, are very hard to handle mathematically. The
landmark work of Christodoulou and Klainerman on small data global existence
for the vacuum Einstein equations \cite{christodoulou} is so complicated as to
discourage any attempts to incorporate matter. The more recent alternative
proof of Lindblad and Rodnianski \cite{lindblad} looks much more promising.
Nevertheless, it seems to require good decay estimates for higher derivatives,
i.e. estimates similar to those proved for the Vlasov-Poisson system here.

The Vlasov-Poisson system in $N$ dimensions reads
\begin{align}
f_{t}+v\cdot\nabla_{x}f+\gamma\nabla_{x}\phi\cdot\nabla_{v}f  &
=0\;\;,\;\;\;x\in\mathbb{R}^{N}\;\;,t>0\label{VP1}\\
\Delta\phi &  =\int_{\mathbb{R}^{N}}fdv\equiv\rho\left(  x,t\right)
\;\;,\;\;x\in\mathbb{R}^{N}\;\;,t>0 \label{VP2}%
\end{align}
where $f=f\left(  x,v,t\right)  .\;$ In the following we assume that $f\left(
x,v,0\right)  =f_{0}\left(  x,v\right)  $ has finite $L^{1}$ norm
and\textbf{\ }$N\geq3$\textbf{.} The sign $\gamma=\pm1$ corresponds to the
plasma physics and gravitational problem respectively. Since the results of
this paper apply equally to both cases, we will restrict our analysis to the
case $\gamma=1.$ No sign condition on $f_{0}$ is needed. However, some
additional decay properties for $f_{0}\left(  x,v\right)  $ will be assumed.

Global existence and decay estimates for the Vlasov-Poisson system were
studied in \cite{bardos} in three spatial dimensions under suitable smallness
and regularity assumptions for $f_{0}.$ These estimates are optimal in the
rate of decay for the density $\rho$ since, for small compactly supported
initial data, the volume of the support of $\rho$ can be bounded by $C\left(
1+t\right)  ^{3}$ so that if the decay in $L^{\infty}$ was stronger than
$\left(  1+t\right)  ^{-3},$ the total number of particles (i.e. the $L^{1}$ 
norm of $\rho$) would decay, leading to a contradiction due to the 
conservation of that quantity. However, they do not provide the optimal rate 
of decay for the derivatives that could be expected on dimensional grounds.

For small initial data the dynamics of the Vlasov-Poisson system might be
expected to be dominated by the free streaming part of the equation:
\[
f_{t}+v\cdot\nabla_{x}f=0
\]
because the term $\nabla_{x}\phi\cdot\nabla_{v}f$ is quadratic in the density.
(Actually this is a consequence of the Bardos-Degond analysis). If we assume 
that the dynamics of the problem is dominated by the free streaming regime as
$t\rightarrow\infty$ and the initial density of particles is, say, compactly
supported (fast enough decay works similarly), the velocities of the particles
would be bounded by a number of order one. Therefore, the support of the
density $\rho$ would spread linearly. The field $\nabla\phi$ generated by a
particle density with finite mass spread over a region of order $t$ decreases
as $\frac{1}{t^{2}}$ as can be easily seen by means of a rescaling argument.
Notice that a posteriori this provides a justification for the assumption
that was made before concerning the finiteness of the deviation of the
velocities of the particles due to the interaction of the field.

The main contribution of this paper is the development of a technique that
allows us to obtain optimal decay estimates for the solutions of the VP
system\textbf{\ }in $N$-dimensional space. More precisely, the rescaling
argument sketched above suggests that the particles spread into a region of
volume $t^{N}$ in the $x$-coordinate. Since the total mass of the particles is
of order one it would be natural to expect the following estimates for the
density
\begin{align*}
\left\vert \rho\right\vert  &  \leq\frac{C}{\left(  t+1\right)  ^{N}}\\
\left\vert \nabla\rho\right\vert  &  \leq\frac{C}{\left(  t+1\right)  ^{N+1}%
}\\
\left\vert \nabla^{2}\rho\right\vert  &  \leq\frac{C}{\left(  t+1\right)
^{N+2}}\\
&  ...\\
\left\vert \nabla^{k}\rho\right\vert  &  \leq\frac{C}{\left(  t+1\right)
^{N+k}}%
\end{align*}

The first estimate was obtained by Bardos-Degond for the case $N=3$ and can be
similarly extended to the case $N>3$. Our method allows us to obtain the
corresponding estimates for the derivatives for small initial data.

The basic idea of the method is as follows. It is easy to see self-similar
behaviour for the density (and the derivatives) in the free streaming case.
Indeed, in that case, integration along characteristics yields
\[
f\left(  x,v,t\right)  =f_{0}\left(  x-vt,v\right)
\]
whence:
\[
\rho\left(  x,t\right)  =\int f\left(  x,v,t\right)  dv=\int f_{0}\left(
x-vt,v\right)  dv
\]

In order to obtain self-similar behaviour we make the change of variables
\begin{align*}
x_{0}  &  =x-vt\\
dv  &  =\left\vert \det\left(  \frac{\partial v}{\partial x_{0}}\right)
\right\vert dx_{0}=\frac{1}{t^{N}}dx_{0}%
\end{align*}
whence:
\[
\rho\left(  x,t\right)  =\frac{1}{t^{N}}\int f_{0}\left(  x_{0},\frac{x-x_{0}%
}{t}\right)  dx_{0}%
\]
In the limit $t\rightarrow\infty$ this formula yields the self-similar
behaviour in the region where $\left\vert x\right\vert $ is of order $t$:
\begin{equation}
\rho\left(  x,t\right)  \sim\frac{1}{t^{N}}\int f_{0}\left(  x_{0},\frac{x}%
{t}\right)  dx_{0}=\frac{1}{t^{N}}\rho_{fs}\left(  \frac{x}{t}\right)
\label{S0E1}%
\end{equation}
Here the asymptotic free streaming density $\rho_{fs}$ is given by
\[
\rho_{fs}\left(  y\right)  \equiv\int f_{0}\left(  x_{0},y\right)  dx_{0}%
\]

Notice that (\ref{S0E1}), at least formally, provides the desired estimates
for the derivatives of $\rho.$ The key idea of our argument is a method for
generalizing this method to the full VP system with small initial data.
The main point is the following. Suppose that the
characteristics starting at $x_{0},$ $v_{0}$ reach the points $x,$ $v$ at time
$t.$ Assuming suitable invertibility conditions any pair of variables in the
set $\left(  x_{0},v_{0},x,v\right)  $ can be used as a set of independent
variables in order to represent the others. The previous argument for the free
streaming case suggests using $x,\;x_{0}$ as independent variables. However,
in order to determine the functions that provide $v_{0},$ $v$ in terms of
$x,\;x_{0}$ it turns out to be necessary to solve a boundary value problem
for the characteristic equations. The main argument of this paper consists in
proving that such a boundary problem can be solved for small initial densities
and that the corresponding solutions of such a boundary value problem satisfy
suitable regularity and decay estimates.

Using a similar method it is possible to obtain not only estimates for the
derivatives of the density, but also convergence of the solutions of the VP 
system to a self-similar solution. More precisely, we rewrite the problem 
using the self-similar variables $y=\frac{x}{\left(  t+1\right)  },\;v=v,\;\tau
=\log\left(  t+1\right)  ,\;f=\frac{1}{\left(  t+1\right)  ^{N}}g$ and after
integrating the resulting equations along characteristics we replace the
variables $\left(  y,v\right)  $ by $\left(  y,y_{0}\right)  .$ This change of
variables requires the solution of a boundary value problem for the 
characteristic
equations analogous to the one described above. Such a boundary value problem
can be analyzed in detail as $\tau\rightarrow\infty,$ and this provides the
asymptotic behaviour of the density $\rho$ as $t\rightarrow\infty.$ One of the
relevant results of the analysis is the fact that although the asymptotics of
the solutions is self-similar, the precise function describing the asymptotics
of the density depends in a very sensitive manner on the choice of the initial
data $f_{0}\left(  x,v\right)  .$ This analysis\ is done in the last section
of the paper.

The paper is organized as follows. In Section 2, we derive estimates for the
density and its derivatives. In Section 3, we prove convergence to the
self-similar solution. Throughout the paper, $C>0$ will denote a generic 
constant that may change from line to line and is independent of 
$t,$ $\varepsilon_{0}$, $f_{0}$.

\section{ESTIMATING $\rho$ AND ITS DERIVATIVES.}

\subsection{The main result.}

We will use  the following function spaces extensively\textbf{.}%

\begin{align}
&  \left\Vert \rho\right\Vert _{X_{k,\alpha}}\label{Xka}\\
&  =\sup_{t\geq0}\{\int\left\vert \rho\left(  x,t\right)  \right\vert
d^{N}x+\left(  t+1\right)  ^{N}\sum_{\ell=0}^{k}\left(  t+1\right)  ^{\ell
}\left\Vert \nabla^{\ell}\rho\left(  \cdot,t\right)  \right\Vert _{L^{\infty
}\left(  \mathbb{R}^{N}\right)  }\nonumber\\
&  +\left(  t+1\right)  ^{N+k+\alpha}\sup_{x,y\in\mathbb{R}^{N}}%
\frac{\left\vert \nabla^{k}\rho\left(  x,t\right)  -\nabla^{k}\rho\left(
y,t\right)  \right\vert }{\left\vert x-y\right\vert ^{\alpha}}\}\nonumber
\end{align}
where $0<\alpha<1.$

\bigskip%

\begin{align}
\left\|  \phi\right\|  _{Y_{k,\alpha}}  &  =\sup_{t\geq0}\{\left(  t+1\right)
^{N-2}\sum_{\ell=1}^{k+2}\left(  t+1\right)  ^{\ell}\left\|  \nabla^{\ell}%
\phi\left(  \cdot,t\right)  \right\|  _{L^{\infty}\left(  \mathbb{R}%
^{N}\right)  }\label{Yka}\\
&  +\left(  t+1\right)  ^{N-2+k+\alpha}\sup_{x,x^{\prime}\in\mathbb{R}^{N}%
}\frac{\left|  \nabla^{k+2}\phi\left(  x,t\right)  -\nabla^{k+2}\phi\left(
x^{\prime},t\right)  \right|  }{\left|  x-x^{\prime}\right|  ^{\alpha}%
}\}\nonumber
\end{align}

The main result of the paper is the following Theorem:

\begin{theorem}
\label{thm}Suppose that $f_{0}\left(  x,v\right)  $ satisfies the following
assumptions:
\begin{align}
\sum_{\ell=0}^{k}\sum_{m=0}^{\ell}\left\vert \frac{\partial^{\ell}f_{0}%
}{\partial x^{m}\partial v^{\ell-m}}\right\vert  &  \leq\frac{\varepsilon_{0}%
}{\left(  1+\left\vert x\right\vert \right)  ^{K}\left(  1+\left\vert
v\right\vert \right)  ^{K}}\;\;,\;\;\label{smalldata}\\
\sum_{m=0}^{k}\sup_{x,x^{\prime}\in\mathbb{R}^{N}}\frac{\left\vert
\frac{\partial^{k}f_{0}}{\partial x^{m}\partial v^{k-m}}\left(  x,v\right)
-\frac{\partial^{k}f_{0}}{\partial x^{m}\partial v^{k-m}}\left(  x^{\prime
},v\right)  \right\vert }{\left\vert x-x^{\prime}\right\vert ^{\alpha}}  &
\leq\frac{\varepsilon_{0}}{\left(  1+\left\vert v\right\vert \right)  ^{K}%
}\;,\;0<\alpha<1\;,\;\nonumber\\
\sum_{m=0}^{k}\sup_{\left\vert v^{\prime}-v\right\vert \leq1}\frac{\left\vert
\frac{\partial^{k}f_{0}}{\partial x^{m}\partial v^{k-m}}\left(  x,v\right)
-\frac{\partial^{k}f_{0}}{\partial x^{m}\partial v^{k-m}}\left(  x,v^{\prime
}\right)  \right\vert }{\left\vert v-v^{\prime}\right\vert ^{\alpha}}  &
\leq\frac{\varepsilon_{0}}{\left(  1+\left\vert x\right\vert \right)
^{K}\left(  1+\left\vert v\right\vert \right)  ^{K}}\;,\;0<\alpha
<1\;,\;\nonumber
\end{align}
for some suitable $K>N$ and $\varepsilon_{0}>0$ small enough.
Then there exists a corresponding solution of the Vlasov-Poisson system with
\[
\left\|  \rho\right\|  _{X_{k,\alpha}}\leq C\varepsilon_{0}%
\]

\end{theorem}

\bigskip

A result analogous to Theorem \ref{thm} was proved in the case $k=0$ under
slightly different assumptions on $f_{0},$ by Bardos-Degond (cf.
\cite{bardos}). The main contribution of this paper is to derive the optimal
decay estimates for the derivatives of $\rho$. \bigskip

\subsection{A basic boundary value problem for the characteristic curves.}

We introduce some basic notation\textbf{. }Suppose that the characteristics
starting at $\left(  x_{0},v_{0}\right)  $ at time $t=0$ reach the point
$\left(  x,v\right)  $ at time $t.$ The basic idea in the paper is to use
$x,\;x_{0}$ as independent variables to describe the values of $v$\ and
$v_{0}$. More precisely, we will write:
\begin{align}
v  &  =w\left(  x,x_{0},t\right) \label{cv1}\\
v_{0}  &  =w_{0}\left(  x,x_{0},t\right)  \label{cv2}%
\end{align}
Note that the existence of the functions\textbf{\ }$w,$ $w_{0}$ is assured by
the implicit function theorem as well as the estimates that we will derive
later where it will be shown that the change of variables is close to a change
of variables that can be inverted explicitly. By changing the variable
$v$ to $x_{0}$ in the integral with $x$ and $t$ fixed, it then follows, using
(\ref{cv1}) that
\[
dv=\left|  \det\left(  \frac{\partial w\left(  x,x_{0},t\right)  }{\partial
x_{0}}\right)  \right|  dx_{0}%
\]
whence
\[
\rho\left(  x,t\right)  =\int_{\mathbb{R}^{N}}f\left(  x,v,t\right)  =\int
f\left(  x_{0},v_{0}\right)  dv=\int f\left(  x_{0},w_{0}\left(
x,x_{0},t\right)  \right)  \left|  \det\frac{\partial w\left(  x,x_{0}%
,t\right)  }{\partial x_{0}}\right|  dx_{0}.
\]

We now formulate the following auxiliary boundary value problem that describes
the evolution of the characteristics starting at the spatial point $x_{0}$ at
the initial time and reaching the point $x$ at time $t:$%
\begin{equation}
\frac{dX\left(  s\right)  }{ds}=V\left(  s\right)  ,\ \frac{dV\left(
s\right)  }{ds}=\nabla\phi\left(  X\left(  s\right)  ,s\right)  ,\ X\left(
t\right)  =x,\ X\left(  0\right)  =x_{0} \label{XV1}%
\end{equation}
Notice that the functions\textbf{\ }$X\left(  s\right)  ,\;V\left(  s\right)
$\textbf{\ }depend also on the variables\textbf{\ }$x,\;x_{0},\;t.$ However,
for simplicity, we will not write the dependence on these variables explicitly
unless it is needed. We then rewrite the above characteristics as a
perturbation from those associated to the free streaming case as follows
\begin{equation}
\frac{dX\left(  s\right)  }{ds}=V\left(  s\right)  =\frac{x-x_{0}}{t}%
+\varphi\left(  s\right)  ,\;\;\ \frac{d\varphi}{ds}=\nabla\phi\left(
X\left(  s\right)  ,s\right)  ,\ X\left(  t\right)  =x,\ X\left(  0\right)
=x_{0}, \label{XV}%
\end{equation}
where $\varphi\left(  s\right)  =\varphi\left(  s;x,x_{0},t\right)  $ is the
perturbed value of the velocity with respect to the free streaming case.
Notice that in the limit of zero density $\rho\equiv0,$ the field $\phi$
vanishes and $\varphi\left(  s\right)  \equiv0.$

We examine the derivatives of the function $\varphi\left(  s\right)  $ with
respect to the variables\textbf{\ }$x,$\textbf{\ }$x_{0},$ in order to derive
suitable estimates for $\rho.$ The density function $\rho$ can be represented
as
\begin{equation}
\rho\left(  x,t\right)  =\int f\left(  x,v,t\right)  dv=\int f_{0}\left(
x_{0},V\left(  0;x,w\left(  t,x,x_{0}\right)  ,t\right)  \right)  \left\vert
\det\left(  \frac{\partial w}{\partial x_{0}}\right)  \right\vert dx_{0}.
\label{rho}%
\end{equation}
Along the characteristics, we have
\[
\frac{\partial w}{\partial x_{0}}=\frac{\partial V}{\partial x_{0}}\left(
t\right)  =-\frac{1}{t}I_{N}+\frac{\partial\varphi}{\partial x_{0}}\left(
t\right)  ,
\]
where $I_{N}$ is the $N$-dimensional identity matrix.

On the other hand, we wish to obtain estimates for the derivatives of $\rho.$
Suppose for the moment that we restrict our attention to the first derivative
of $\rho$ with respect to $x$. Such a derivative is given by:%

\begin{align*}
\frac{\partial\rho}{\partial x}\left(  x,t\right)   &  =\int\frac{\partial
f_{0}}{\partial v}\left(  x_{0},V\left(  0;x,w\left(  t,x,x_{0}\right)
,t\right)  \right)  \frac{\partial V}{\partial x}\left(  0\right)  \left|
\det\left(  \frac{\partial w}{\partial x_{0}}\right)  \right|  dx_{0}\\
&  +\int f_{0}\left(  x_{0},V\left(  0;x,w\left(  t,x,x_{0}\right)  ,t\right)
\right)  \frac{\partial}{\partial x}\left[  \left|  \det\left(  \frac{\partial
w}{\partial x_{0}}\right)  \right|  \right]  dx_{0}.
\end{align*}
To estimate the first derivative of $\rho$ reduces to derive estimates for:
\[
\frac{\partial V}{\partial x}\left(  s=0\right)  ,\ \frac{\partial V}{\partial
x_{0}}\left(  s=t\right)  ,\ \frac{\partial^{2}V}{\partial x\partial x_{0}%
}\left(  s=t\right)  .
\]
Equivalently
\[
\frac{\partial\varphi}{\partial x}\left(  0\right)  ,\ \frac{\partial\varphi
}{\partial x_{0}}\left(  t\right)  ,\ \frac{\partial^{2}\varphi}{\partial
x\partial x_{0}}\left(  t\right)  .
\]

Notice that the equation of the characteristics (\ref{XV}) indicates that in
order to obtain bounds for two derivatives with respect to $x,$\ $x_{0}$ of
the characteristic curves we need to estimate three derivatives of the
potential $\phi.$ These are the exact number of derivatives that can be
expected to be estimated from the Poisson equation under the assumption that
$\frac{\partial\rho}{\partial x}$ is bounded. Nevertheless, in order to avoid
the standard problems that arise in the regularity estimates for the Poisson
equation in the spaces $C^{k},$ it is necessary to work with the H\"{o}lder
spaces $C^{k,\alpha}$.

\subsection{Estimates on the regularity and the rate of decay of
$\varphi\left(  s;x,x_{0},t\right)  $ in terms of the properties of the
potential $\phi.$}

We present a key a priori estimate for $\varphi\,\ $in terms of $\phi$ in the
following Proposition. We define two norms with respect to the spatial
variable $x$.

\begin{definition}
For $u\left(  \cdot\right)  \in L^{\infty}\left(  \mathbb{R}^{N}\right)  ,$%
\[
\left\|  u\right\|  _{L^{\infty}\left(  x\right)  }\equiv\operatorname*{ess}%
\sup_{x\in\mathbb{R}^{N}}\left|  u\left(  x\right)  \right|  ,
\]
\ 

For $u\left(  \cdot\right)  \in C^{0,\alpha}\left(  \mathbb{R}^{N}\right)  ,$%
\[
\left[  u\right]  _{0,\alpha,\left(  x\right)  }\equiv\sup_{x_{1},x_{2}%
\in\mathbb{R}^{N}}\frac{\left|  u\left(  x_{1}\right)  -u\left(  x_{2}\right)
\right|  }{\left|  x_{1}-x_{2}\right|  ^{\alpha}}.
\]

\end{definition}

For notational simplicity, in the following, we will use $\left\Vert u\left(
s\right)  \right\Vert _{L^{\infty}\left(  x\right)  },\left[  u\left(
s\right)  \right]  _{0,\alpha,\left(  x\right)  }$ instead of $\left\Vert
u\left(  s;\cdot,x_{0},t\right)  \right\Vert _{L^{\infty}\left(  x\right)
},\left[  u\left(  s;\cdot,x_{0},t\right)  \right]  _{0,\alpha,\left(
x\right)  }$, which in fact depend on $s,x_{0},$ and $t.$ For example, $\left[
\frac{\partial\varphi}{\partial x}\left(  s\right)  \right]  _{0,\alpha
,\left(  x\right)  }$ will denote $\left[  \frac{\partial\varphi}{\partial
x}\left(  s;\cdot,x_{0},t\right)  \right]  _{0,\alpha,\left(  x\right)  }.$

\begin{proposition}
\label{regul}Suppose that
\[
\left\|  \phi\right\|  _{Y_{1,\alpha}}\leq\varepsilon_{0}.
\]
for a suitable $\varepsilon_{0}>0$ sufficiently small. Suppose that $t\geq1.$
Then, the following a priori estimate holds:
\begin{align*}
t\sup_{0\leq s\leq t}\left\|  \frac{\partial\varphi}{\partial x}\left(
s\right)  \right\|  _{L^{\infty}\left(  x\right)  }+t^{1+\alpha}\sup_{0\leq
s\leq t}\left[  \frac{\partial\varphi}{\partial x}\left(  s\right)  \right]
_{0,\alpha,\left(  x\right)  }  &  \leq C\left\|  \phi\right\|  _{Y_{1,\alpha
}}\\
t\left\|  \frac{\partial\varphi}{\partial x_{0}}\left(  t\right)  \right\|
_{L^{\infty}\left(  x\right)  }+t^{1+\alpha}\left[  \frac{\partial\varphi
}{\partial x_{0}}\left(  t\right)  \right]  _{0,\alpha,\left(  x\right)
}+t^{2}\left\|  \frac{\partial^{2}\varphi}{\partial x\partial x_{0}}\left(
t\right)  \right\|  _{L^{\infty}\left(  x\right)  }+t^{2+\alpha}\left[
\frac{\partial^{2}\varphi}{\partial x\partial x_{0}}\left(  t\right)  \right]
_{0,\alpha,\left(  x\right)  }  &  \leq C\left\|  \phi\right\|  _{Y_{1,\alpha
}}%
\end{align*}%
\begin{align*}
&  \int_{0}^{t}\left\|  \frac{\partial\varphi}{\partial x_{0}}\left(
s\right)  \right\|  _{L^{\infty}\left(  x\right)  }ds+t^{\alpha}\int_{0}%
^{t}\left[  \frac{\partial\varphi}{\partial x_{0}}\left(  s\right)  \right]
_{0,\alpha,\left(  x\right)  }ds\\
&  +t\int_{0}^{t}ds\left\|  \frac{\partial^{2}\varphi}{\partial x\partial
x_{0}}\left(  s\right)  \right\|  _{L^{\infty}\left(  x\right)  }+t^{1+\alpha
}\int_{0}^{t}\left[  \frac{\partial^{2}\varphi}{\partial x\partial x_{0}%
}\left(  s\right)  \right]  _{0,\alpha,\left(  x\right)  }ds\\
&  \leq C\left\|  \phi\right\|  _{Y_{1,\alpha}}%
\end{align*}
for some suitable constant $C>0$ independent of $t,\;\varepsilon_{0}.$
\end{proposition}

Proposition \ref{regul} is the main new technical result of the paper. This 
estimate
provides optimal decay properties for the derivatives of the characteristics
in terms of the decay properties of the derivatives of the potential $\phi.$

Furthermore, we derive a generalization of Proposition \ref{regul} under
additional regularity and decay assumptions for the potential $\phi:$

\begin{proposition}
\label{regulhigh}Suppose that
\[
\left\|  \phi\right\|  _{Y_{\ell,\alpha}}\leq\varepsilon_{0},\;\;\ell\geq2
\]
for a suitable $\varepsilon_{0}>0$ sufficiently small. Suppose that $t\geq1.$
Then, the following estimates hold:
\begin{align}
\sum_{k=1}^{\ell}t^{k}\sup_{0\leq s\leq t}\left\|  \frac{\partial^{k}\varphi
}{\partial x^{k}}\left(  s\right)  \right\|  _{L^{\infty}\left(  x\right)
}+t^{\ell+\alpha}\sup_{0\leq s\leq t}\left[  \frac{\partial^{\ell}\varphi
}{\partial x^{\ell}}\left(  s\right)  \right]  _{0,\alpha,\left(  x\right)  }
&  \leq C\left\|  \phi\right\|  _{Y_{\ell,\alpha}}\label{phixl-inf}\\
\sum_{k=1}^{\ell}t^{k}\int_{0}^{t}\left\|  \frac{\partial^{k+1}\varphi
}{\partial x^{k}\partial x_{0}}\left(  s\right)  \right\|  _{L^{\infty}\left(
x\right)  }ds+t^{\ell+\alpha}\int_{0}^{t}\left[  \frac{\partial^{\ell
+1}\varphi}{\partial x^{\ell}\partial x_{0}}\left(  s\right)  \right]
_{0,\alpha,\left(  x\right)  }ds  &  \leq C\left\|  \phi\right\|
_{Y_{\ell,\alpha}}\label{phixx0l}\\
\sum_{k=1}^{\ell}t^{k+1}\left\|  \frac{\partial^{k+1}\varphi}{\partial
x^{k}\partial x_{0}}\left(  s\right)  \right\|  _{L^{\infty}\left(  x\right)
}+t^{\ell+1+\alpha}\left[  \frac{\partial^{\ell+1}\varphi}{\partial x^{\ell
}\partial x_{0}}\left(  t\right)  \right]  _{0,\alpha,\left(  x\right)  }  &
\leq C\left\|  \phi\right\|  _{Y_{\ell,\alpha}} \label{phixx0lt}%
\end{align}
for some constant $C>0$ independent of $t,\;\varepsilon_{0}.$
\end{proposition}

\bigskip\bigskip

\subsubsection{Preliminary results: Integral equation satisfied by
$\varphi\left(  s,x,x_{0},t\right)  .$}

The perturbed velocity $\varphi\left(  s\right)  $ satisfies the 
following integral equation

\begin{lemma}
$\varphi\left(  s\right)  =\varphi\left(  s;x,x_{0},t\right)  $ in (\ref{XV})
satisfies the integral equation
\begin{equation}
\varphi\left(  s\right)  =-\int_{s}^{t}G\left(  \xi\right)  d\xi+\frac{1}%
{t}\int_{0}^{t}\xi G\left(  \xi\right)  d\xi, \label{phi}%
\end{equation}
where
\begin{align}
G\left(  \xi\right)   &  \equiv\nabla\phi\left(  X\left(  \xi\right)
,\xi\right)  ,\ \label{G}\\
X\left(  \xi\right)   &  =x_{0}+\frac{x-x_{0}}{t}\xi+\int_{0}^{\xi}%
\varphi\left(  \bar{s}\right)  d\bar{s}. \label{X}%
\end{align}

\end{lemma}

\begin{proof}
We integrate (\ref{XV}) to obtain
\begin{equation}
\varphi\left(  s\right)  =\varphi\left(  0\right)  +\int_{0}^{s}G\left(
\xi\right)  d\xi,\ \;\;X\left(  s\right)  =x_{0}+\frac{x-x_{0}}{t}s+\int
_{0}^{s}\varphi\left(  \bar{s}\right)  d\bar{s}. \label{phiX}%
\end{equation}
The boundary condition $X\left(  t\right)  =x$ yields
\begin{equation}
\int_{0}^{t}\varphi\left(  s\right)  ds=0. \label{zerointegral}%
\end{equation}
Therefore, integrating the first equation in (\ref{phiX}) from $s=0$ to $s=t$
and using (\ref{zerointegral}), we have
\[
\varphi\left(  0\right)  =-\frac{1}{t}\int_{0}^{t}\left[  \int_{0}^{s}G\left(
\xi\right)  d\xi\right]  ds.
\]
Thus, we obtain from (\ref{phiX})
\[
\varphi\left(  s\right)  =-\frac{1}{t}\int_{0}^{t}\left[  \int_{0}^{s}G\left(
\xi\right)  d\xi\right]  ds+\int_{0}^{s}G\left(  \xi\right)  d\xi.
\]
By changing the order of integration in the first integral on the right-hand
side of this formula, we deduce (\ref{phi}).
\end{proof}

\bigskip

Note that
\begin{align*}
V\left(  t\right)   &  =V\left(  t;x,x_{0},t\right)  =v\equiv w\left(
x,x_{0},t\right)  =\frac{x-x_{0}}{t}+\varphi\left(  t;x,x_{0},t\right)  ,\\
V\left(  0\right)   &  =V\left(  0;x,x_{0},t\right)  =v_{0}\equiv w_{0}\left(
x,x_{0},t\right)  =\frac{x-x_{0}}{t}+\varphi\left(  0;x,x_{0},t\right)  .
\end{align*}

The integral equation (\ref{phi}), (\ref{G}), (\ref{X}) is the key ingredient
that will be used to derive optimal regularity and decay estimates for the
functions $w\left(  x,x_{0},t\right)  ,\;V_{0}\left(  x,x_{0},t\right)  .$ As
a first step we prove that the solutions of (\ref{phi}), (\ref{G}), (\ref{X})
are well defined if $\varepsilon_{0}$ is small enough.

\begin{lemma}
\label{Solvability}(Solvability) Let $k\geq0$ be an integer. There exists
$\varepsilon_{0}>0$ such that, for any $t\geq1,\;x\in\mathbb{R}^{N},\;x_{0}%
\in\mathbb{R}^{N}$ and any function $\phi$ satisfying
\[
\left\|  \phi\right\|  _{Y_{k,\alpha}}\leq\varepsilon_{0},
\]
there exists a unique solution $\varphi\left(  \cdot\right)  =\varphi\left(
\cdot;x,x_{0},t\right)  \in C(\left[  0,t\right]  )$ of (\ref{phi}),
(\ref{G}), (\ref{X}).
\end{lemma}

\begin{proof}
Let the space of functions for $\varphi\left(  s\right)  $ be
\[
\mathcal{X}\equiv\left\{  \varphi\in C(\left[  0,t\right]  ):\sup_{0\leq s\leq
t}\left\vert \varphi\left(  s\right)  \right\vert \leq1\right\}  .
\]
Let
\begin{equation}
\mathcal{J}\left(  \varphi\right)  \left(  s\right)  \equiv-\int_{s}%
^{t}G_{\varphi}\left(  \xi\right)  d\xi+\frac{1}{t}\int_{0}^{t}\xi G_{\varphi
}\left(  \xi\right)  d\xi, \label{Jdef}%
\end{equation}
where
\[
G_{\varphi}\left(  \xi\right)  =\nabla\phi\left(  x_{0}+\frac{x-x_{0}}{t}%
\xi+\int_{0}^{\xi}\varphi\left(  \bar{s}\right)  d\bar{s},\xi\right)
\]
We first show that $\mathcal{J}$ is a well-defined operator in the
space\textbf{\ }$\mathcal{X}$. By definition, we have
\[
\left\vert G_{\varphi}\left(  \xi\right)  \right\vert \leq\left\Vert
\nabla\phi\left(  \cdot,\xi\right)  \right\Vert _{L^{\infty}\left(  x\right)
}\leq\frac{\left\Vert \phi\right\Vert _{Y_{k,\alpha}}}{\left(  \xi+1\right)
^{N-1}}.
\]
This yields
\begin{align*}
\left\vert \mathcal{J}\left(  \varphi\right)  \left(  s\right)  \right\vert
&  \leq\int_{s}^{t}\left\vert G_{\varphi}\left(  \xi\right)  \right\vert
d\xi+\frac{1}{t}\int_{0}^{t}\xi\left\vert G_{\varphi}\left(  \xi\right)
\right\vert d\xi\\
&  \leq\left\Vert \phi\right\Vert _{Y_{k,\alpha}}\left[  \int_{0}^{t}%
\frac{d\xi}{\left(  1+\xi\right)  ^{N-1}}+\frac{1}{t}\int_{0}^{t}\frac{\xi
}{\left(  1+\xi\right)  ^{N-1}}d\xi\right] \\
&  \leq\left\Vert \phi\right\Vert _{Y_{k,\alpha}}\left[  1+\frac{\log\left(
1+t\right)  }{t}\right]  \leq C\left\Vert \phi\right\Vert _{Y_{k,\alpha}},
\end{align*}
where we have estimated the last integral term using the fact that
$N\geq3.$ This idea will be used repeatedly in the following. \textbf{\ }Thus
$\mathcal{J}\left(  \varphi\right)  $ is bounded for all $t$ and $\mathcal{J}$
is well-defined in the space\textbf{\ }$\mathcal{X}$. We next show that
\thinspace the operator $\mathcal{J}$ is contractive. Taking the difference of
$\mathcal{J}\left(  \varphi_{1}\right)  $ and $\mathcal{J}\left(  \varphi
_{1}\right)  $ yields
\[
\lbrack\mathcal{J}\left(  \varphi_{1}\right)  -\mathcal{J}\left(  \varphi
_{2}\right)  ]\left(  s\right)  =-\int_{s}^{t}\left[  G_{\varphi_{1}}\left(
\xi\right)  -G_{\varphi_{2}}\left(  \xi\right)  \right]  d\xi+\frac{1}{t}%
\int_{0}^{t}\xi\left[  G_{\varphi_{1}}\left(  \xi\right)  -G_{\varphi_{2}%
}\left(  \xi\right)  \right]  d\xi.
\]
Using the definition of $\left\Vert \phi\right\Vert _{Y_{k,\alpha}},$ we have
\begin{align*}
&  \sup_{0\leq s\leq t}\left\vert [\mathcal{J}\left(  \varphi_{1}\right)
-\mathcal{J}\left(  \varphi_{2}\right)  ]\left(  s\right)  \right\vert \\
&  \leq\int_{s}^{t}\left\Vert \nabla^{2}\phi\left(  \cdot,\xi\right)
\right\Vert _{L^{\infty}\left(  x\right)  }\left[  \int_{0}^{\xi}\left\vert
\varphi_{1}\left(  \bar{s}\right)  -\varphi_{2}\left(  \bar{s}\right)
\right\vert d\bar{s}\right]  d\xi+\frac{1}{t}\int_{0}^{t}\xi\left\Vert
\nabla^{2}\phi\left(  \cdot,\xi\right)  \right\Vert _{L^{\infty}\left(
x\right)  }\left[  \int_{0}^{\xi}\left\vert \varphi_{1}\left(  \bar{s}\right)
-\varphi_{2}\left(  \bar{s}\right)  \right\vert d\bar{s}\right]  d\xi\\
&  \leq\left\Vert \phi\right\Vert _{Y_{k,\alpha}}\left[  \int_{s}^{t}\frac
{1}{\left(  \xi+1\right)  ^{N}}\left[  \int_{0}^{\xi}\left\vert \varphi
_{1}\left(  \bar{s}\right)  -\varphi_{2}\left(  \bar{s}\right)  \right\vert
d\bar{s}\right]  d\xi+\frac{1}{t}\int_{0}^{t}\frac{\xi}{\left(  \xi+1\right)
^{N}}\left[  \int_{0}^{\xi}\left\vert \varphi_{1}\left(  \bar{s}\right)
-\varphi_{2}\left(  \bar{s}\right)  \right\vert d\bar{s}\right]  d\xi\right]
\\
&  \leq C\varepsilon_{0}\sup_{0\leq s\leq t}\left\vert \left(  \varphi
_{1}-\varphi_{2}\right)  \left(  s\right)  \right\vert +C\varepsilon_{0}%
\frac{\log\left(  t+1\right)  }{t}\sup_{0\leq s\leq t}\left\vert \left(
\varphi_{1}-\varphi_{2}\right)  \left(  s\right)  \right\vert \\
&  \leq C\varepsilon_{0}\sup_{0\leq s\leq t}\left\vert \left(  \varphi
_{1}-\varphi_{2}\right)  \left(  s\right)  \right\vert .
\end{align*}
Thus we can choose $\varepsilon_{0}$ small enough such that $C\varepsilon
_{0}<1$ in the above inequality and conclude that $\mathcal{J}$ is
contractive. Notice that $C$ is independent of $t,\;x_{0},x$ and
$\varepsilon_{0}$ can be chosen independently of these variables. Then by the
Banach fixed point theorem, we deduce the existence and uniqueness for 
$\varphi$ satisfying (\ref{phi}) in the space\textbf{\ }$\mathcal{X}$.
\end{proof}

\begin{corollary}
Let $\varphi\left(  s\right)  $ be the solution in Lemma \ref{Solvability}.
Then we have
\[
\int_{0}^{t}\varphi\left(  s\right)  ds=0.
\]

\end{corollary}

\begin{proof}
Using (\ref{Jdef}) we derive the identity:
\[
\int_{0}^{t}\mathcal{J}\left(  \varphi\right)  ds=0.
\]
Therefore, using Lemma \ref{Solvability}, we have
\[
\varphi=\mathcal{J}\left(  \varphi\right)  .
\]
This completes the proof of Corollary.
\end{proof}

Notice that this Corollary implies that the characteristics (\ref{XV}) 
satisfy the
desired boundary condition for $X\left(  s\right)  ,$ i.e., $X\left(
t\right)  =x$.

We now turn to decay estimates for the density function $\rho$ and its
derivatives. Before we proceed, we state some basic properties of the
H\"{o}lder norms.

\begin{lemma}
\label{Schauder}
\[
\left[  fg\right]  _{0,\alpha,\left(  x\right)  }\leq C\{\left\|  f\right\|
_{L^{\infty}}[g]_{0,\alpha,\left(  x\right)  }+\left\|  g\right\|
_{L^{\infty}}[f]_{0,\alpha,\left(  x\right)  }\}\;,
\]
for any $f,g\in L^{\infty}\cap C^{0,\alpha},$%
\[
\lbrack f]_{0,\alpha,\left(  x\right)  }\leq C\left\|  f\right\|  _{L^{\infty
}}^{1-\alpha}\left\|  \nabla f\right\|  _{L^{\infty}}^{\alpha},
\]
for any $f\in W^{1,\infty},$%
\[
\left[  F\circ u\right]  _{0,\alpha,\left(  x\right)  }\leq\left[  F\right]
_{0,\alpha,\left(  x\right)  }\left\|  \nabla u\right\|  _{L^{\infty}}%
^{\alpha},
\]
for any $F\in C^{0,\alpha}$ and any $u\in W^{1,\infty}.$

\begin{proof}
The results in the lemma are standard estimates for H\"{o}lder norms.
\end{proof}
\end{lemma}

\bigskip

\subsubsection{The proof of Proposition \ref{regul}.}

The proof of Proposition \ref{regul} follows from a sequence of lemmas. There
are three ideas that will appear repeatedly in all the remaining arguments of
this paper. Estimating terms like $\frac{\partial\varphi}{\partial x_{0}}$ and
its derivatives with respect to $x,$ it is not possible to obtain bounds for
the rate of decay suggested by dimensional considerations for all the values
of $s\in\left[  0,t\right]  .$ It is, however,  possible to obtain such
optimal decay estimates for the integrals of such terms in the interval
$\left[  0,t\right]  $ as well as for the time $s=t$ which is the only one
where such optimal estimates are really needed. The second idea is that it is
convenient to obtain, before deriving pointwise estimates, integral estimates
for terms like $\int_{0}^{t}\left\|  \frac{\partial\varphi}{\partial x_{0}%
}\left(  s\right)  \right\|  _{L^{\infty}\left(  x\right)  }ds$. The third
idea is that the estimates for terms that do not contain derivatives with
respect $x_{0}$ are more easily obtained by directly estimating the supremum
over the interval $\left[  0,t\right]  $ and using Gronwall-type arguments,
without any need for estimating integrals over the interval $\left[
0,t\right]  $.

\begin{lemma}
\label{phidecayk0}\bigskip There exists $\varepsilon_{0}$ small such that for
$t>1$ and any function $\phi$ satisfying
\[
\left\|  \phi\right\|  _{Y_{0,\alpha}}\leq\varepsilon_{0},
\]
we have
\begin{equation}
\int_{0}^{t}\left\|  \frac{\partial\varphi}{\partial x_{0}}\left(  s\right)
\right\|  _{L^{\infty}\left(  x\right)  }ds\leq C\left\|  \phi\right\|
_{Y_{0,\alpha}}, \label{phix01}%
\end{equation}%
\begin{equation}
\left\|  \frac{\partial\varphi}{\partial x_{0}}\left(  t\right)  \right\|
_{L^{\infty}\left(  x\right)  }\leq C\frac{\left\|  \phi\right\|
_{Y_{0,\alpha}}}{t}. \label{phix02}%
\end{equation}

\end{lemma}

\begin{proof}
Differentiating (\ref{phi}) with respect to $x_{0}$ yields
\begin{equation}
\frac{\partial\varphi}{\partial x_{0}}\left(  s\right)  =-\int_{s}^{t}%
\frac{\partial}{\partial x_{0}}G\left(  \xi\right)  d\xi+\frac{1}{t}\int
_{0}^{t}\xi\frac{\partial}{\partial x_{0}}G\left(  \xi\right)  d\xi.
\label{phi_x0}%
\end{equation}
where, for simplicity\textbf{\ }$G_{\varphi}=G.$ We now take $\frac{\partial
}{\partial x_{0}}$ of (\ref{G})-(\ref{X}) to get
\[
\frac{\partial}{\partial x_{0}}G\left(  \xi\right)  =\nabla^{2}\phi\left(
X\left(  \xi\right)  ,\xi\right)  \frac{\partial}{\partial x_{0}}X\left(
\xi\right)  =\nabla^{2}\phi\left(  X\left(  \xi\right)  ,\xi\right)
[(1-\frac{\xi}{t})I+\int_{0}^{\xi}\frac{\partial\phi}{\partial x_{0}}]
\]
and since $\frac{\xi}{t}\leq1,$ we have
\begin{align*}
\left\|  \frac{\partial}{\partial x_{0}}G\left(  \xi\right)  \right\|
_{L^{\infty}\left(  x\right)  }  &  \leq\frac{\left\|  \phi\right\|
_{Y_{0,\alpha}}}{\left(  \xi+1\right)  ^{N}}\left\|  (1-\frac{\xi}{t}%
)I+\int_{0}^{\xi}\frac{\partial\phi}{\partial x_{0}}\left(  \bar{s}\right)
d\bar{s}\right\|  _{L^{\infty}\left(  x\right)  }\\
&  \leq\frac{\left\|  \phi\right\|  _{Y_{0,\alpha}}}{\left(  \xi+1\right)
^{N}}\left[  1+\int_{0}^{\xi}\left\|  \frac{\partial\varphi}{\partial x_{0}%
}\left(  \bar{s}\right)  \right\|  _{L^{\infty}\left(  x\right)  }d\bar
{s}\right]  .
\end{align*}
Putting the above into (\ref{phi_x0}) yields
\begin{align}
\left\|  \frac{\partial\varphi}{\partial x_{0}}\left(  s\right)  \right\|
_{L^{\infty}\left(  x\right)  }  &  \leq C\left\|  \phi\right\|
_{Y_{0,\alpha}}\int_{s}^{t}\frac{1}{\left(  \xi+1\right)  ^{N}}[1+\int
_{0}^{\xi}\left\|  \frac{\partial\varphi}{\partial x_{0}}\left(  \bar
{s}\right)  \right\|  _{L^{\infty}\left(  x\right)  }d\bar{s}]d\xi
\label{phi_x01}\\
&  +C\left\|  \phi\right\|  _{Y_{0,\alpha}}\frac{1}{t}\int_{0}^{t}\frac{\xi
}{\left(  \xi+1\right)  ^{N}}[1+\int_{0}^{\xi}\left\|  \frac{\partial\varphi
}{\partial x_{0}}\left(  \bar{s}\right)  \right\|  _{L^{\infty}\left(
x\right)  }d\bar{s}]d\xi.\nonumber
\end{align}
By integrating the above from $0$ to $t$ and by the assumption, we obtain
\begin{align*}
\int_{0}^{t}\left\|  \frac{\partial\varphi}{\partial x_{0}}\left(  s\right)
\right\|  _{L^{\infty}\left(  x\right)  }ds  &  \leq C\left\|  \phi\right\|
_{Y_{0,\alpha}}\int_{0}^{t}\left[  \int_{s}^{t}\frac{1}{\left(  \xi+1\right)
^{N}}[1+\int_{0}^{\xi}\left\|  \frac{\partial\varphi}{\partial x_{0}}\left(
\bar{s}\right)  \right\|  _{L^{\infty}\left(  x\right)  }d\bar{s}]d\xi\right]
ds\\
&  +C\left\|  \phi\right\|  _{Y_{0,\alpha}}\int_{0}^{t}\frac{\xi}{\left(
\xi+1\right)  ^{N}}[1+\int_{0}^{\xi}\left\|  \frac{\partial\varphi}{\partial
x_{0}}\left(  \bar{s}\right)  \right\|  _{L^{\infty}\left(  x\right)  }%
d\bar{s}]d\xi\\
&  \leq C\left\|  \phi\right\|  _{Y_{0,\alpha}}+C\left\|  \phi\right\|
_{Y_{0,\alpha}}\int_{0}^{t}\left\|  \frac{\partial\varphi}{\partial x_{0}%
}\left(  s\right)  \right\|  _{L^{\infty}\left(  x\right)  }ds\\
&  \leq C\left\|  \phi\right\|  _{Y_{0,\alpha}}+C\varepsilon_{0}\int_{0}%
^{t}\left\|  \frac{\partial\varphi}{\partial x_{0}}\left(  s\right)  \right\|
_{L^{\infty}\left(  x\right)  }ds,
\end{align*}
where we have used the estimate $\int_{0}^{\xi}\left\|  \frac{\partial\varphi
}{\partial x_{0}}\left(  \bar{s}\right)  \right\|  _{L^{\infty}\left(
x\right)  }d\bar{s}\leq\int_{0}^{t}\left\|  \frac{\partial\varphi}{\partial
x_{0}}\left(  \bar{s}\right)  \right\|  _{L^{\infty}\left(  x\right)  }%
d\bar{s}$ and the fact that $\int_{0}^{t}\left[  \int_{s}^{t}\frac{d\xi
}{\left(  \xi+1\right)  ^{3}}\right]  ds=\int_{0}^{t}\frac{1}{\left(
\xi+1\right)  ^{3}}\left[  \int_{0}^{\xi}ds\right]  d\xi\leq C.$ Thus if
$\varepsilon_{0}$ is small enough so that $C\varepsilon_{0}\leq1/2,$ we get
\begin{equation}
\int_{0}^{t}\left\|  \frac{\partial\varphi}{\partial x_{0}}\left(  s\right)
\right\|  _{L^{\infty}\left(  x\right)  }ds\leq C\left\|  \phi\right\|
_{Y_{0,\alpha}}. \label{phi_x02}%
\end{equation}
Putting $s=t$ in (\ref{phi_x01}) and using (\ref{phi_x02}) yields
\begin{align*}
\left\|  \frac{\partial\varphi}{\partial x_{0}}\left(  t\right)  \right\|
_{L^{\infty}\left(  x\right)  }  &  \leq C\left\|  \phi\right\|
_{Y_{0,\alpha}}\frac{1}{t}\int_{0}^{t}\frac{\xi}{\left(  \xi+1\right)  ^{N}%
}[1+\int_{0}^{\xi}\left\|  \frac{\partial\varphi}{\partial x_{0}}\left(
\bar{s}\right)  \right\|  _{L^{\infty}\left(  x\right)  }d\bar{s}]d\xi\\
&  \leq C\frac{\left\|  \phi\right\|  _{Y_{0,\alpha}}}{t}(1+\left\|
\phi\right\|  _{Y_{0,\alpha}})\leq C\frac{\left\|  \phi\right\|
_{Y_{0,\alpha}}}{t}.
\end{align*}
Thus we obtain (\ref{phix01}) and (\ref{phix02}). This completes the proof of
the lemma.
\end{proof}

\begin{lemma}
\label{phidecayk1}There exists $\varepsilon_{0}$ small such that for $\left\|
\phi\right\|  _{Y_{1,\alpha}}\leq\varepsilon_{0}$ and\textbf{\ }$t>1$ we have
the following decay estimates
\begin{equation}
\sup_{0\leq s\leq t}\left\|  \frac{\partial\varphi}{\partial x}\left(
s\right)  \right\|  _{L^{\infty}\left(  x\right)  }\leq C\frac{\left\|
\phi\right\|  _{Y_{1,\alpha}}}{t}\;\;,\ \label{phix1}%
\end{equation}%
\begin{equation}
\int_{0}^{t}\left\|  \frac{\partial^{2}\varphi}{\partial x\partial x_{0}%
}\left(  s\right)  \right\|  _{L^{\infty}\left(  x\right)  }ds\leq
C\frac{\left\|  \phi\right\|  _{Y_{1,\alpha}}}{t}\;\;,\ \label{phixx0}%
\end{equation}%
\begin{equation}
\left\|  \frac{\partial^{2}\varphi}{\partial x\partial x_{0}}\left(  t\right)
\right\|  _{L^{\infty}\left(  x\right)  }\leq C\frac{\left\|  \phi\right\|
_{Y_{1,\alpha}}}{t^{2}}\;\;. \label{phixx01}%
\end{equation}

\end{lemma}

\begin{proof}
Differentiating (\ref{G})-(\ref{X}) with respect to $x$ we get
\[
\frac{\partial}{\partial x}G\left(  \xi\right)  =\nabla^{2}\phi\left(
X\left(  \xi\right)  ,\xi\right)  \frac{\partial}{\partial x}X\left(
\xi\right)  =\nabla^{2}\phi\left(  X\left(  \xi\right)  ,\xi\right)
[\frac{\xi}{t}I+\int_{0}^{\xi}\frac{\partial\phi}{\partial x}\left(  \bar
{s}\right)  d\bar{s}]
\]
and thus
\begin{align*}
\left\|  \frac{\partial}{\partial x}G\left(  \xi\right)  \right\|
_{L^{\infty}\left(  x\right)  }  &  \leq\frac{\left\|  \phi\right\|
_{Y_{1,\alpha}}}{\left(  \xi+1\right)  ^{N}}\left\|  \frac{\xi}{t}I+\int
_{0}^{\xi}\frac{\partial\phi}{\partial x}\left(  \bar{s}\right)  d\bar
{s}\right\|  _{L^{\infty}\left(  x\right)  }\\
&  \leq\frac{\left\|  \phi\right\|  _{Y_{1,\alpha}}}{\left(  \xi+1\right)
^{N-1}}\left[  \frac{\xi}{\left(  \xi+1\right)  t}+\frac{1}{\left(
\xi+1\right)  }\int_{0}^{\xi}\left\|  \frac{\partial\varphi}{\partial
x}\left(  \bar{s}\right)  \right\|  _{L^{\infty}\left(  x\right)  }d\bar
{s}\right] \\
&  \leq\frac{\left\|  \phi\right\|  _{Y_{1,\alpha}}}{\left(  \xi+1\right)
^{N-1}}\left[  \frac{1}{t}+\frac{1}{\left(  \xi+1\right)  }\int_{0}^{\xi
}\left\|  \frac{\partial\varphi}{\partial x}\left(  \bar{s}\right)  \right\|
_{L^{\infty}\left(  x\right)  }d\bar{s}\right]  .
\end{align*}
Differentiating (\ref{phi}) with respect to $x$ and using the estimate above,
we obtain
\begin{align*}
\left\|  \frac{\partial\varphi}{\partial x}\left(  s\right)  \right\|
_{L^{\infty}\left(  x\right)  }  &  \leq C\left\|  \phi\right\|
_{Y_{1,\alpha}}\int_{s}^{t}\frac{1}{\left(  \xi+1\right)  ^{N-1}}\left[
\frac{1}{t}+\frac{1}{\left(  \xi+1\right)  }\int_{0}^{\xi}\left\|
\frac{\partial\varphi}{\partial x}\left(  \bar{s}\right)  \right\|
_{L^{\infty}\left(  x\right)  }d\bar{s}\right]  d\xi\\
&  +C\left\|  \phi\right\|  _{Y_{1,\alpha}}\frac{1}{t}\int_{0}^{t}\frac{\xi
}{\left(  \xi+1\right)  ^{N-1}}\left[  \frac{1}{t}+\frac{1}{\left(
\xi+1\right)  }\int_{0}^{\xi}\left\|  \frac{\partial\varphi}{\partial
x}\left(  \bar{s}\right)  \right\|  _{L^{\infty}\left(  x\right)  }d\bar
{s}\right]  d\xi\\
&  \leq C\frac{\left\|  \phi\right\|  _{Y_{1,\alpha}}}{t}+C\left\|
\phi\right\|  _{Y_{1,\alpha}}\sup_{0\leq s\leq t}\left\|  \frac{\partial
\varphi}{\partial x}\left(  s\right)  \right\|  _{L^{\infty}\left(  x\right)
}\\
&  +C\left\|  \phi\right\|  _{Y_{1,\alpha}}\frac{\log\left(  t+1\right)
}{t^{2}}+C\left\|  \phi\right\|  _{Y_{1,\alpha}}\frac{\log\left(  t+1\right)
}{t}\sup_{0\leq s\leq t}\left\|  \frac{\partial\varphi}{\partial x}\left(
s\right)  \right\|  _{L^{\infty}\left(  x\right)  }.
\end{align*}
It then follows that, for\textbf{\ }$t>1:$
\begin{equation}
\sup_{0\leq s\leq t}\left\|  \frac{\partial\varphi}{\partial x}\left(
s\right)  \right\|  _{L^{\infty}\left(  x\right)  }\leq C\frac{\left\|
\phi\right\|  _{Y_{1,\alpha}}}{t}+C\left\|  \phi\right\|  _{Y_{1,\alpha}}%
\sup_{0\leq s\leq t}\left\|  \frac{\partial\varphi}{\partial x}\left(
s\right)  \right\|  _{L^{\infty}\left(  x\right)  } \label{phix}%
\end{equation}
By the assumption, if $\varepsilon_{0}$ is small enough, then we get
\[
\sup_{0\leq s\leq t}\left\|  \frac{\partial\varphi}{\partial x}\left(
s\right)  \right\|  _{L^{\infty}\left(  x\right)  }\leq C\frac{\left\|
\phi\right\|  _{Y_{1,\alpha}}}{t}.
\]
and (\ref{phix1}) follows.

In order to derive (\ref{phixx0}) and (\ref{phixx01}), we compute $\frac
{\partial^{2}G}{\partial x\partial x_{0}},\;\frac{\partial^{2}\varphi
}{\partial x\partial x_{0}}$ using (\ref{G}), (\ref{X}) and (\ref{phiX}):
\begin{align}
\frac{\partial^{2}}{\partial x\partial x_{0}}G\left(  \xi\right)   &
=\nabla^{3}\phi\left(  X\left(  \xi\right)  ,\xi\right)  \frac{\partial
}{\partial x}X\left(  \xi\right)  \frac{\partial}{\partial x_{0}}X\left(
\xi\right)  +\nabla^{2}\phi\left(  X\left(  \xi\right)  ,\xi\right)
\frac{\partial^{2}}{\partial x\partial x_{0}}X\left(  \xi\right) \nonumber\\
&  =\nabla^{3}\phi\left(  X\left(  \xi\right)  ,\xi\right)  [\frac{\xi}%
{t}I+\int_{0}^{\xi}\frac{\partial\phi}{\partial x}\left(  \bar{s}\right)
d\bar{s}][(1-\frac{\xi}{t})I+\int_{0}^{\xi}\frac{\partial\phi}{\partial x_{0}%
}\left(  \bar{s}\right)  d\bar{s}]\nonumber\\
&  +\nabla^{2}\phi\left(  X\left(  \xi\right)  ,\xi\right)  \int_{0}^{\xi
}\frac{\partial^{2}\varphi}{\partial x\partial x_{0}}\left(  \bar{s}\right)
d\bar{s}, \label{dersG}%
\end{align}

\begin{equation}
\frac{\partial^{2}\varphi}{\partial x\partial x_{0}}\left(  s\right)
=-\int_{s}^{t}\frac{\partial^{2}G}{\partial x\partial x_{0}}\left(
\xi\right)  d\xi+\frac{1}{t}\int_{0}^{t}\xi\frac{\partial^{2}G}{\partial
x\partial x_{0}}\left(  \xi\right)  d\xi. \label{dersephi}%
\end{equation}
Taking the norm $\left\|  \cdot\right\|  _{L^{\infty}\left(  x\right)  }$ of
this equation, integrating the resulting formula with respect to $s,$ using
Lemma \ref{phidecayk0}, (\ref{phix}), and the definition of $\left\|
\phi\right\|  _{Y_{1,\alpha}},$ we obtain
\begin{align*}
&  \int_{0}^{t}\left\|  \frac{\partial^{2}\varphi}{\partial x\partial x_{0}%
}\left(  s\right)  \right\|  _{L^{\infty}\left(  x\right)  }ds\\
&  \leq\int_{0}^{t}\left[  \int_{s}^{t}\left\|  \frac{\partial^{2}G}{\partial
x\partial x_{0}}\left(  \xi\right)  \right\|  _{L^{\infty}\left(  x\right)
}d\xi\right]  ds+\int_{0}^{t}\xi\left\|  \frac{\partial^{2}G}{\partial
x\partial x_{0}}\left(  \xi\right)  \right\|  _{L^{\infty}\left(  x\right)
}d\xi\ \\
&  \leq C\frac{\left\|  \phi\right\|  _{Y_{1,\alpha}}}{t}\int_{0}^{t}%
ds\int_{s}^{t}\frac{\xi d\xi}{\left(  \xi+1\right)  ^{N+3}}+C\left\|
\phi\right\|  _{Y_{1,\alpha}}\left(  \int_{0}^{t}\left\|  \frac{\ \partial
^{2}\varphi}{\partial x\partial x_{0}}\left(  s\right)  \right\|  _{L^{\infty
}\left(  x\right)  }ds\right)  \left(  \int_{0}^{t}\left[  \int_{s}^{t}%
\frac{d\xi}{\left(  \xi+1\right)  ^{N}}\right]  ds\right) \\
&  +C\frac{\left\|  \phi\right\|  _{Y_{1,\alpha}}}{t}\int_{0}^{t}\frac
{d\xi\ \xi^{2}}{\left(  \xi+1\right)  ^{N+1}}+C\left\|  \phi\right\|
_{Y_{1,\alpha}}\left(  \int_{0}^{t}\left\|  \frac{\ \partial^{2}\varphi
}{\partial x\partial x_{0}}\left(  s\right)  \right\|  _{L^{\infty}\left(
x\right)  }ds\right)  \left(  \int_{0}^{t}\frac{\xi d\xi}{\left(
\xi+1\right)  ^{N}}\right) \\
&  \leq C\frac{\left\|  \phi\right\|  _{Y_{1,\alpha}}}{t}+C\left\|
\phi\right\|  _{Y_{1,\alpha}}\int_{0}^{t}\left\|  \frac{\partial^{2}\varphi
}{\partial x\partial x_{0}}\left(  s\right)  \right\|  _{L^{\infty}\left(
x\right)  }ds.
\end{align*}
Thus if $\varepsilon_{0}$ is small enough, we get
\[
\int_{0}^{t}\left\|  \frac{\partial^{2}\varphi}{\partial x\partial x_{0}%
}\left(  s\right)  \right\|  _{L^{\infty}\left(  x\right)  }ds\leq
C\frac{\left\|  \phi\right\|  _{Y_{1,\alpha}}}{t}
\]
and (\ref{phixx0}) follows. We now set $s=t$ in (\ref{dersephi}) and use
(\ref{dersG}) to obtain
\begin{align*}
\left\|  \frac{\partial^{2}\varphi}{\partial x\partial x_{0}}\left(  t\right)
\right\|  _{L^{\infty}\left(  x\right)  }  &  \leq\frac{1}{t}\int_{0}^{t}%
\xi\left\|  \frac{\partial^{2}G}{\partial x\partial x_{0}}\left(  \xi\right)
\right\|  _{L^{\infty}\left(  x\right)  }d\xi\\
&  \leq\frac{C\left\|  \phi\right\|  _{Y_{1,\alpha}}}{t^{2}}\int_{0}^{t}%
\frac{\xi^{2}d\xi}{\left(  \xi+1\right)  ^{N+1}}+\frac{C\left\|  \phi\right\|
_{Y_{1,\alpha}}}{t}\int_{0}^{t}\frac{\xi d\xi}{\left(  \xi+1\right)  ^{N}}%
\int_{0}^{\xi}\left\|  \frac{\partial^{2}\varphi}{\partial x\partial x_{0}%
}\left(  \bar{s}\right)  \right\|  _{L^{\infty}\left(  x\right)  }d\bar{s}\\
&  \leq\frac{C\left\|  \phi\right\|  _{Y_{1,\alpha}}}{t^{2}}+\frac{C\left\|
\phi\right\|  _{Y_{1,\alpha}}^{2}}{t^{2}}\leq C\frac{\left\|  \phi\right\|
_{Y_{1,\alpha}}}{t^{2}}\;\;
\end{align*}
and the proof of Lemma \ref{phidecayk1} is complete.
\end{proof}

\bigskip

In order to complete the proof of Proposition \ref{regul}, it only remains to
obtain estimates for the H\"{o}lder seminorms of $\frac{\partial\varphi
}{\partial x},$ $\frac{\partial\varphi}{\partial x_{0}},\,\frac{\partial
^{2}\varphi}{\partial x\partial x_{0}}.$ These bounds are obtained using ideas
analogous to those used in the two previous lemmas.

\begin{lemma}
There exists $\varepsilon_{0}$ small such that for $t>1$ and $\left\|
\phi\right\|  _{Y_{1,\alpha}}\leq\varepsilon_{0},$ we have the following decay
estimates:
\begin{equation}
\sup_{0\leq s\leq t}\left[  \frac{\partial\varphi}{\partial x}\left(
s\right)  \right]  _{0,\alpha,\left(  x\right)  }\leq C\frac{\left\|
\phi\right\|  _{Y_{1,\alpha}}}{t^{1+\alpha}}, \label{phixa}%
\end{equation}%
\begin{equation}
\int_{0}^{t}\left[  \frac{\partial\varphi}{\partial x_{0}}\left(  s\right)
\right]  _{0,\alpha,\left(  x\right)  }ds\leq C\frac{\left\|  \phi\right\|
_{Y_{1,\alpha}}}{t^{\alpha}}, \label{phix0-int}%
\end{equation}%
\begin{equation}
\left[  \frac{\partial\varphi}{\partial x_{0}}\left(  t\right)  \right]
_{0,\alpha,\left(  x\right)  }\leq C\frac{\left\|  \phi\right\|
_{Y_{1,\alpha}}}{t^{1+\alpha}}, \label{phix03}%
\end{equation}%
\[
\int_{0}^{t}\left[  \frac{\partial^{2}\varphi}{\partial x\partial x_{0}%
}\left(  s\right)  \right]  _{0,\alpha,\left(  x\right)  }ds\leq
C\frac{\left\|  \phi\right\|  _{Y_{1,\alpha}}}{t^{1+\alpha}},
\]%
\begin{equation}
\ \left[  \frac{\partial^{2}\varphi}{\partial x\partial x_{0}}\left(
t\right)  \right]  _{0,\alpha,\left(  x\right)  }\leq C\frac{\left\|
\phi\right\|  _{Y_{1,\alpha}}}{t^{2+\alpha}}. \label{phixx0a}%
\end{equation}
\textbf{\ }
\end{lemma}

\begin{proof}
Using Lemma \ref{Schauder} and (\ref{phix1}), we get
\begin{align*}
\left[  \frac{\partial}{\partial x}G\left(  \xi\right)  \right]
_{0,\alpha,\left(  x\right)  }  &  \leq C\left\|  \nabla^{2}\phi\left(
X\left(  \xi\right)  ,\xi\right)  \right\|  _{L^{\infty}\left(  x\right)
}\left[  \frac{\xi}{t}I+\int_{0}^{\xi}\frac{\partial\phi}{\partial x}\left(
\bar{s}\right)  d\bar{s}\right]  _{0,\alpha,\left(  x\right)  }\\
&  +C\left[  \nabla^{2}\phi\right]  _{0,\alpha,\left(  x\right)  }\left\|
\frac{\xi}{t}I+\int_{0}^{\xi}\frac{\partial\phi}{\partial x}\left(  \bar
{s}\right)  d\bar{s}\right\|  _{L^{\infty}\left(  x\right)  }^{\alpha}\left\|
\frac{\xi}{t}I+\int_{0}^{\xi}\frac{\partial\phi}{\partial x}\right\|
_{L^{\infty}\left(  x\right)  }\\
&  \leq\frac{C\left\|  \phi\right\|  _{Y_{1,\alpha}}}{\left(  \xi+1\right)
^{N}}\int_{0}^{\xi}\left[  \frac{\partial\phi}{\partial x}\left(  \bar
{s}\right)  \right]  _{0,\alpha,\left(  x\right)  }d\bar{s}+\frac{C\left\|
\phi\right\|  _{Y_{1,\alpha}}}{\left(  \xi+1\right)  ^{N+\alpha}}\frac
{\xi^{\alpha}}{t^{\alpha}}\frac{\xi}{t}\\
&  \leq\frac{C\left\|  \phi\right\|  _{Y_{1,\alpha}}}{\left(  \xi+1\right)
^{N}}\int_{0}^{\xi}\left[  \frac{\partial\phi}{\partial x}\left(  \bar
{s}\right)  \right]  _{0,\alpha,\left(  x\right)  }d\bar{s}+\frac{C\left\|
\phi\right\|  _{Y_{1,\alpha}}}{\left(  \xi+1\right)  ^{N-1}t^{1+\alpha}}.
\end{align*}
Differentiating (\ref{phi}) with respect to $x$, taking the H\"{o}lder norm,
and using the previous estimate, we get
\begin{align*}
\left[  \frac{\partial\varphi}{\partial x}\left(  s\right)  \right]
_{0,\alpha,\left(  x\right)  }  &  \leq\int_{s}^{t}\left[  \frac{\partial
}{\partial x}G\left(  \xi\right)  \right]  _{0,\alpha,\left(  x\right)  }%
d\xi+\frac{1}{t}\int_{0}^{t}\xi\left[  \frac{\partial}{\partial x}G\left(
\xi\right)  \right]  _{0,\alpha,\left(  x\right)  }d\xi\\
&  \leq\frac{C\left\|  \phi\right\|  _{Y_{1,\alpha}}}{t^{1+\alpha}}\int
_{s}^{t}\frac{d\xi}{\left(  \xi+1\right)  ^{N-1}}+C\left\|  \phi\right\|
_{Y_{1,\alpha}}\sup_{0\leq s\leq t}\left[  \frac{\partial\phi}{\partial
x}\left(  s\right)  \right]  _{0,\alpha,\left(  x\right)  }\int_{s}^{t}%
\frac{\xi d\xi}{\left(  \xi+1\right)  ^{N}}\\
&  +\frac{C\left\|  \phi\right\|  _{Y_{1,\alpha}}}{t^{2+\alpha}}\int_{0}%
^{t}\frac{\xi d\xi}{\left(  \xi+1\right)  ^{N-1}}+\frac{C\left\|
\phi\right\|  _{Y_{1,\alpha}}}{t}\sup_{0\leq s\leq t}\left[  \frac
{\partial\phi}{\partial x}\left(  s\right)  \right]  _{0,\alpha,\left(
x\right)  }\int_{s}^{t}\frac{\xi^{2}d\xi}{\left(  \xi+1\right)  ^{N}}\\
&  \leq\frac{C\left\|  \phi\right\|  _{Y_{1,\alpha}}}{t^{1+\alpha}}+C\left\|
\phi\right\|  _{Y_{1,\alpha}}\sup_{0\leq s\leq t}\left[  \frac{\partial\phi
}{\partial x}\left(  s\right)  \right]  _{0,\alpha,\left(  x\right)  },
\end{align*}
where we have used that $\left[  \frac{\xi}{t}I+\int_{0}^{\xi}\frac
{\partial\phi}{\partial x}\left(  \bar{s}\right)  d\bar{s}\right]
_{0,\alpha,\left(  x\right)  }=\left[  \int_{0}^{\xi}\frac{\partial\phi
}{\partial x}\left(  \bar{s}\right)  d\bar{s}\right]  _{0,\alpha,\left(
x\right)  }$ as well as the fact that, due to (\ref{phix1}), $\left\|
\frac{\xi}{t}I+\int_{0}^{\xi}\frac{\partial\phi}{\partial x}\right\|
_{L^{\infty}\left(  x\right)  }\leq\frac{2\xi}{t}.$ We then deduce
(\ref{phixa}) provided $\varepsilon_{0}$ is small enough. We now derive the
H\"{o}lder estimate of $\frac{\partial\varphi}{\partial x_{0}}\left(
t\right)  $. By interpolation, the H\"{o}lder inequality, Lemma 
\ref{phidecayk0}, and Lemma \ref{phidecayk1}, we obtain the following two 
estimates
\begin{align*}
&  \int_{0}^{t}\left[  \frac{\partial\varphi}{\partial x_{0}}\left(  s\right)
\right]  _{0,\alpha,\left(  x\right)  }ds\\
&  \leq C\int_{0}^{t}\left\|  \frac{\partial\varphi}{\partial x_{0}}\left(
s\right)  \right\|  _{L^{\infty}\left(  x\right)  }^{1-\alpha}\left\|
\frac{\partial^{2}\varphi}{\partial x_{0}\partial x}\left(  s\right)
\right\|  _{L^{\infty}\left(  x\right)  }^{\alpha}ds\\
&  \leq C\left(  \int_{0}^{t}\left\|  \frac{\partial\varphi}{\partial x_{0}%
}\left(  s\right)  \right\|  _{L^{\infty}\left(  x\right)  }ds\right)
^{1-\alpha}\left(  \int_{0}^{t}\left\|  \frac{\partial^{2}\varphi}{\partial
x_{0}\partial x}\left(  s\right)  \right\|  _{L^{\infty}\left(  x\right)
}ds\right)  ^{\alpha}\\
&  \leq C\left\|  \phi\right\|  _{Y_{1,\alpha}}^{1-\alpha}\frac{\left\|
\phi\right\|  _{Y_{1,\alpha}}^{\alpha}}{t^{\alpha}}\leq C\frac{\left\|
\phi\right\|  _{Y_{1,\alpha}}}{t^{\alpha}}.
\end{align*}%
\[
\left[  \frac{\partial\varphi}{\partial x_{0}}\left(  t\right)  \right]
_{0,\alpha,\left(  x\right)  }\leq C\left\|  \frac{\partial\varphi}{\partial
x_{0}}\left(  t\right)  \right\|  _{L^{\infty}\left(  x\right)  }^{1-\alpha
}\left\|  \frac{\partial^{2}\varphi}{\partial x_{0}\partial x}\left(
t\right)  \right\|  _{L^{\infty}\left(  x\right)  }^{\alpha}\leq
C\frac{\left\|  \phi\right\|  _{Y_{1,\alpha}}^{1-\alpha}}{t^{1-\alpha}}%
\frac{\left\|  \phi\right\|  _{Y_{1,\alpha}}^{\alpha}}{t^{2\alpha}}\leq
C\frac{\left\|  \phi\right\|  _{Y_{1,\alpha}}}{t^{1+\alpha}}.
\]
Therefore, (\ref{phix0-int}) and (\ref{phix03}) follow.

We now use Lemma \ref{Schauder}, Lemma \ref{phidecayk0}, Lemma
\ref{phidecayk1}, and (\ref{phixa})-(\ref{phix03}) to get
\begin{align*}
&  \left[  \frac{\partial^{2}}{\partial x\partial x_{0}}G\left(  \xi\right)
\right]  _{0,\alpha,\left(  x\right)  }\\
&  \leq C\left\Vert \nabla^{3}\phi\left(  X\left(  \xi\right)  ,\xi\right)
\right\Vert _{L^{\infty}\left(  x\right)  }\left[  (1-\frac{\xi}{t})I+\int
_{0}^{\xi}\frac{\partial\phi}{\partial x_{0}}\right]  _{0,\alpha,\left(
x\right)  }\left\Vert \frac{\xi}{t}I+\int_{0}^{\xi}\frac{\partial\phi
}{\partial x}\right\Vert _{L^{\infty}\left(  x\right)  }\\
&  +C\left\Vert \nabla^{3}\phi\left(  X\left(  \xi\right)  ,\xi\right)
\right\Vert _{L^{\infty}\left(  x\right)  }\left\Vert (1-\frac{\xi}{t}%
)I+\int_{0}^{\xi}\frac{\partial\phi}{\partial x_{0}}\right\Vert _{L^{\infty
}\left(  x\right)  }\left[  \frac{\xi}{t}I+\int_{0}^{\xi}\frac{\partial\phi
}{\partial x}\right]  _{0,\alpha,\left(  x\right)  }\\
&  +C\left[  \nabla^{3}\phi\right]  _{0,\alpha,\left(  x\right)  }\left\Vert
\frac{\xi}{t}I+\int_{0}^{\xi}\frac{\partial\phi}{\partial x}\right\Vert
_{L^{\infty}\left(  x\right)  }^{\alpha}\left\Vert (1-\frac{\xi}{t})I+\int
_{0}^{\xi}\frac{\partial\phi}{\partial x_{0}}\right\Vert _{L^{\infty}\left(
x\right)  }\left\Vert \frac{\xi}{t}I+\int_{0}^{\xi}\frac{\partial\phi
}{\partial x}\right\Vert _{L^{\infty}\left(  x\right)  }\\
&  +C\left[  \nabla^{2}\phi\right]  _{0,\alpha,\left(  x\right)  }\left\Vert
\frac{\xi}{t}I+\int_{0}^{\xi}\frac{\partial\phi}{\partial x}\right\Vert
_{L^{\infty}\left(  x\right)  }^{\alpha}\int_{0}^{\xi}\left\Vert
\frac{\partial^{2}\varphi}{\partial x\partial x_{0}}\left(  \bar{s}\right)
\right\Vert _{L^{\infty}\left(  x\right)  }d\bar{s}\\
&  +C\left\Vert \nabla^{2}\phi\left(  X\left(  \xi\right)  ,\xi\right)
\right\Vert _{L^{\infty}\left(  x\right)  }\int_{0}^{\xi}\left[
\frac{\partial^{2}\varphi}{\partial x\partial x_{0}}\left(  \bar{s}\right)
\right]  _{0,\alpha,\left(  x\right)  }d\bar{s}\\
&  \leq\frac{C\left\Vert \phi\right\Vert _{Y_{1,\alpha}}}{\left(
\xi+1\right)  ^{N+1}}\frac{1}{t^{\alpha}}\frac{\xi}{t}+\frac{C\left\Vert
\phi\right\Vert _{Y_{1,\alpha}}}{\left(  \xi+1\right)  ^{N+1+\alpha}}\frac
{\xi^{\alpha}}{t^{\alpha}}\frac{\xi}{t}+\frac{C\left\Vert \phi\right\Vert
_{Y_{1,\alpha}}}{\left(  \xi+1\right)  ^{N+\alpha}}\frac{\xi^{\alpha}%
}{t^{\alpha}}\frac{C\left\Vert \phi\right\Vert _{Y_{1,\alpha}}}{t}\\
&  +\frac{C\left\Vert \phi\right\Vert _{Y_{1,\alpha}}}{\left(  \xi+1\right)
^{N}}\int_{0}^{\xi}\left[  \frac{\partial^{2}\varphi}{\partial x\partial
x_{0}}\left(  \bar{s}\right)  \right]  _{0,\alpha,\left(  x\right)  }d\bar
{s}\\
&  \leq\frac{C\left\Vert \phi\right\Vert _{Y_{1,\alpha}}}{\left(
\xi+1\right)  ^{N}}\int_{0}^{\xi}\left[  \frac{\partial^{2}\varphi}{\partial
x\partial x_{0}}\left(  s\right)  \right]  _{0,\alpha,\left(  x\right)
}ds+\frac{C\left\Vert \phi\right\Vert _{Y_{1,\alpha}}}{\left(  \xi+1\right)
^{N}t^{1+\alpha}}.
\end{align*}
Similarly$,$ we get from (\ref{phi}) as well as the estimate above
\begin{align}
&  \left[  \frac{\partial^{2}\varphi}{\partial x\partial x_{0}}\left(
s\right)  \right]  _{0,\alpha,\left(  x\right)  }\label{phix0a}\\
&  \leq\int_{s}^{t}\left[  \frac{\partial^{2}}{\partial x\partial x_{0}%
}G\left(  \xi\right)  \right]  _{0,\alpha,\left(  x\right)  }d\xi+\frac{1}%
{t}\int_{0}^{t}\xi\left[  \frac{\partial^{2}}{\partial x\partial x_{0}%
}G\left(  \xi\right)  \right]  _{0,\alpha,\left(  x\right)  }d\xi\nonumber\\
&  \leq\frac{C\left\Vert \phi\right\Vert _{Y_{1,\alpha}}}{t^{1+\alpha}}%
\int_{s}^{t}\frac{d\xi}{\left(  \xi+1\right)  ^{N}}+C\left\Vert \phi
\right\Vert _{Y_{1,\alpha}}\int_{s}^{t}\frac{1}{\left(  \xi+1\right)  ^{N}%
}\left[  \int_{0}^{\xi}\left[  \frac{\partial^{2}\varphi}{\partial x\partial
x_{0}}\left(  \bar{s}\right)  \right]  _{0,\alpha,\left(  x\right)  }d\bar
{s}\right]  d\xi\nonumber\\
&  +\frac{C\left\Vert \phi\right\Vert _{Y_{1,\alpha}}}{t^{2+\alpha}}\int
_{0}^{t}\frac{\xi d\xi}{\left(  \xi+1\right)  ^{N}}+\frac{C\left\Vert
\phi\right\Vert _{Y_{1,\alpha}}}{t}\int_{0}^{t}\frac{\xi}{\left(
\xi+1\right)  ^{N}}\left[  \int_{0}^{\xi}\left[  \frac{\partial^{2}\varphi
}{\partial x\partial x_{0}}\left(  \bar{s}\right)  \right]  _{0,\alpha,\left(
x\right)  }d\bar{s}\right]  d\xi\nonumber
\end{align}
By integrating (\ref{phix0a}) from $s=0$ to $s=t,$ we have
\begin{align*}
\int_{0}^{t}\left[  \frac{\partial^{2}\varphi}{\partial x\partial x_{0}%
}\left(  s\right)  \right]  _{0,\alpha,\left(  x\right)  }ds  &  \leq
\frac{C\left\Vert \phi\right\Vert _{Y_{1,\alpha}}}{t^{1+\alpha}}%
+\frac{C\left\Vert \phi\right\Vert _{Y_{1,\alpha}}}{t^{2+\alpha}}\log\left(
t+1\right)  +C\left\Vert \phi\right\Vert _{Y_{1,\alpha}}\int_{0}^{t}\left[
\frac{\partial^{2}\varphi}{\partial x\partial x_{0}}\left(  s\right)  \right]
_{0,\alpha,\left(  x\right)  }ds\\
&  +\frac{C\left\Vert \phi\right\Vert _{Y_{1,\alpha}}}{t}\log\left(
t+1\right)  \int_{0}^{t}\left[  \frac{\partial^{2}\varphi}{\partial x\partial
x_{0}}\left(  s\right)  \right]  _{0,\alpha,\left(  x\right)  }ds\\
&  \leq\frac{C\left\Vert \phi\right\Vert _{Y_{1,\alpha}}}{t^{1+\alpha}%
}+C\left\Vert \phi\right\Vert _{Y_{1,\alpha}}\int_{0}^{t}\left[
\frac{\partial^{2}\varphi}{\partial x\partial x_{0}}\left(  s\right)  \right]
_{0,\alpha,\left(  x\right)  }ds.
\end{align*}
If $\varepsilon_{0}$ is small enough, then we obtain
\begin{equation}
\int_{0}^{t}\left[  \frac{\partial^{2}\varphi}{\partial x\partial x_{0}%
}\left(  s\right)  \right]  _{0,\alpha,\left(  x\right)  }ds\leq
\frac{C\left\Vert \phi\right\Vert _{Y_{1,\alpha}}}{t^{1+\alpha}}.
\label{phix0a1}%
\end{equation}
Now putting (\ref{phix0a1}) into (\ref{phix0a}) with $s=t$ yields
\[
\left[  \frac{\partial^{2}\varphi}{\partial x\partial x_{0}}\left(  t\right)
\right]  _{0,\alpha,\left(  x\right)  }\leq\frac{C\left\Vert \phi\right\Vert
_{Y_{1,\alpha}}}{t^{2+\alpha}}.
\]
and this completes the proof of the Lemma.
\end{proof}

\subsubsection{The proof of Proposition \ref{regulhigh}.}

Now we prove the decay estimates for the higher order derivatives of $\varphi
$. We prove Proposition \ref{regulhigh} by induction on $\ell$. The induction
hypotheses consist of the following estimates, for $0\leq m<\ell,$
\begin{equation}
\sup_{0\leq s\leq t}\left\|  \frac{\partial^{m}\varphi}{\partial x^{m}}\left(
s\right)  \right\|  _{L^{\infty}\left(  x\right)  }\leq\frac{C\left\|
\phi\right\|  _{Y_{m,\alpha}}}{t^{m}},\ \;\sup_{0\leq s\leq t}\left[
\frac{\partial^{m}\varphi}{\partial x^{m}}\left(  s\right)  \right]
_{0,\alpha,\left(  x\right)  }\leq\frac{C\left\|  \phi\right\|  _{Y_{m,\alpha
}}}{t^{m+\alpha}}, \label{hypo1}%
\end{equation}%
\begin{equation}
\int_{0}^{t}\left\|  \frac{\partial^{m+1}\varphi}{\partial x^{m}\partial
x_{0}}\left(  s\right)  \right\|  _{L^{\infty}\left(  x\right)  }ds\leq
\frac{C\left\|  \phi\right\|  _{Y_{m,\alpha}}}{t^{m}},\ \;\int_{0}^{t}\left[
\frac{\partial^{m+1}\varphi}{\partial x^{m}\partial x_{0}}\left(  s\right)
\right]  _{0,\alpha,\left(  x\right)  }ds\leq\frac{C\left\|  \phi\right\|
_{Y_{m,\alpha}}}{t^{m+\alpha}}, \label{hypo2}%
\end{equation}%
\begin{equation}
\left\|  \frac{\partial^{m+1}\varphi}{\partial x^{m}\partial x_{0}}\left(
t\right)  \right\|  _{L^{\infty}\left(  x\right)  }\leq\frac{C\left\|
\phi\right\|  _{Y_{m,\alpha}}}{t^{m+1}},\ \;\left[  \frac{\partial
^{m+1}\varphi}{\partial x^{m}\partial x_{0}}\left(  t\right)  \right]
_{0,\alpha,\left(  x\right)  }\leq\frac{C\left\|  \phi\right\|  _{Y_{m,\alpha
}}}{t^{m+1+\alpha}}. \label{hypo3}%
\end{equation}
Estimates (\ref{hypo1})-(\ref{hypo3}) have been already proved for
$m=0,1\;$(cf. Proposition \ref{regul}). We begin with the estimates of
$G=G_{\varphi}$ in terms of $\varphi:$

\begin{lemma}
\label{Lem-Gxl}Let $\ell\geq2$ be an integer. Assume the induction hypotheses
(\ref{hypo1})-(\ref{hypo3}). There exists $\varepsilon_{0}$ such that for
$t>1$ and $\left\|  \phi\right\|  _{Y_{\ell,\alpha}}\leq\varepsilon_{0},$ we
have the following
\begin{equation}
\left\|  \frac{\partial^{\ell}G}{\partial x^{\ell}}\left(  \xi\right)
\right\|  _{L^{\infty}\left(  x\right)  }\leq\frac{C\left\|  \phi\right\|
_{Y_{\ell,\alpha}}}{t^{\ell}\left(  \xi+1\right)  ^{N-1}}+\frac{C\left\|
\phi\right\|  _{Y_{\ell,\alpha}}}{\left(  \xi+1\right)  ^{N-1}}\sup_{0\leq
s\leq t}\left\|  \frac{\partial^{\ell}\varphi}{\partial x^{\ell}}\left(
s\right)  \right\|  _{L^{\infty}\left(  x\right)  }\;, \label{Gxl-inf}%
\end{equation}%
\begin{equation}
\left[  \frac{\partial^{\ell}G}{\partial x^{\ell}}\left(  \xi\right)  \right]
_{0,\alpha,\left(  x\right)  }\leq\frac{C\left\|  \phi\right\|  _{Y_{\ell
,\alpha}}}{t^{\ell+\alpha}\left(  \xi+1\right)  ^{N-1}}+\frac{C\left\|
\phi\right\|  _{Y_{\ell,\alpha}}}{\left(  \xi+1\right)  ^{N-1}}\sup_{0\leq
s\leq t}\left[  \frac{\partial^{\ell}\varphi}{\partial x^{\ell}}\left(
s\right)  \right]  _{0,\alpha,\left(  x\right)  }\;, \label{Gxl-a}%
\end{equation}%
\begin{equation}
\left\|  \frac{\partial^{\ell+1}G}{\partial x^{\ell}\partial x_{0}}\left(
\xi\right)  \right\|  _{L^{\infty}\left(  x\right)  }\leq\frac{C\left\|
\phi\right\|  _{Y_{\ell,\alpha}}}{t^{\ell}\left(  \xi+1\right)  ^{N}}%
+\frac{C\left\|  \phi\right\|  _{Y_{\ell,\alpha}}}{\left(  \xi+1\right)  ^{N}%
}\int_{0}^{\xi}\left\|  \frac{\partial^{\ell+1}\varphi}{\partial x^{\ell
}\partial x_{0}}\left(  \bar{s}\right)  \right\|  _{L^{\infty}\left(
x\right)  }d\bar{s} \label{Gxx0l-inf}%
\end{equation}%
\begin{equation}
\left[  \frac{\partial^{\ell+1}G}{\partial x^{\ell}\partial x_{0}}\left(
\xi\right)  \right]  _{0,\alpha,\left(  x\right)  }\leq\frac{C\left\|
\phi\right\|  _{Y_{\ell,\alpha}}}{t^{\ell+\alpha}\left(  \xi+1\right)  ^{N}%
}+\frac{C\left\|  \phi\right\|  _{Y_{\ell,\alpha}}}{\left(  \xi+1\right)
^{N}}\int_{0}^{\xi}\left[  \frac{\partial^{\ell+1}\varphi}{\partial x^{\ell
}\partial x_{0}}\left(  \bar{s}\right)  \right]  _{0,\alpha,\left(  x\right)
}d\bar{s}. \label{Gxx0l-a}%
\end{equation}

\end{lemma}

\begin{proof}
Taking $\frac{\partial^{\ell}}{\partial x^{\ell}}$ of $G$ yields
\begin{align}
\frac{\partial^{\ell}G}{\partial x^{\ell}}\left(  \xi\right)   &
=\sum_{\substack{1\leq i\leq\ell,\\j_{1}+...+j_{i}=\ell}}A_{ij_{1}...j_{i}%
}\frac{\partial^{i}}{\partial x^{i}}\triangledown\phi\left(  X,\xi\right)
\frac{\partial^{j_{1}}}{\partial x^{j_{1}}}X\cdot...\cdot\frac{\partial
^{j_{i}}}{\partial x^{j_{i}}}X\label{Gxl}\\
&  =\triangledown^{2}\phi\left(  X,\xi\right)  \frac{\partial^{\ell}%
X}{\partial x^{\ell}}+\sum_{\substack{j_{m}<\ell\\1\leq m\leq i}%
}...=\triangledown^{2}\phi\left(  X,\xi\right)  \int_{0}^{\xi}\frac
{\partial^{\ell}}{\partial x^{\ell}}\varphi\left(  \bar{s}\right)  d\bar
{s}+\sum_{_{\substack{j_{m}<\ell\\1\leq m\leq i}}}...,\nonumber
\end{align}
where
\[
X=X\left(  \xi\right)  =x_{0}+\frac{x-x_{0}}{t}\xi+\int_{0}^{\xi}%
\varphi\left(  \bar{s}\right)  d\bar{s}.
\]
and where $A_{ij_{1}...j_{i}}$\ are suitable numerical coefficients. Using the
induction hypotheses, we bound each term with $j_{m}<\ell$ for all $1\leq
m\leq i$ on the right-hand side of the above identity as
\[
\left\|  \frac{\partial^{i}}{\partial x^{i}}\triangledown\phi\left(
X,\xi\right)  \frac{\partial^{j_{1}}}{\partial x^{j_{1}}}X\cdot...\cdot
\frac{\partial^{j_{i}}}{\partial x^{j_{i}}}X\right\|  _{L^{\infty}\left(
x\right)  }\leq\frac{C\left\|  \phi\right\|  _{Y_{\ell,\alpha}}}{\left(
\xi+1\right)  ^{i+N-1}}\frac{C\left\|  \phi\right\|  _{Y_{\ell,\alpha}}\xi
}{t^{j_{1}}}\cdot...\cdot\frac{C\left\|  \phi\right\|  _{Y_{\ell,\alpha}}\xi
}{t^{j_{i}}}\leq\frac{C\left\|  \phi\right\|  _{Y_{\ell,\alpha}}}{\left(
\xi+1\right)  ^{N-1}t^{\ell}},
\]%
\begin{align*}
&  \left[  \frac{\partial^{i}}{\partial x^{i}}\triangledown\phi\left(
X,\xi\right)  \frac{\partial^{j_{1}}}{\partial x^{j_{1}}}X\cdot...\cdot
\frac{\partial^{j_{i}}}{\partial x^{j_{i}}}X\right]  _{0,\alpha,\left(
x\right)  }\\
&  \leq\left[  \frac{\partial^{i}}{\partial x^{i}}\triangledown\phi\left(
X,\xi\right)  \right]  _{0,\alpha,\left(  x\right)  }\left\|  \frac
{\partial^{j_{1}}}{\partial x^{j_{1}}}X\cdot...\cdot\frac{\partial^{j_{i}}%
}{\partial x^{j_{i}}}X\right\|  _{L^{\infty}\left(  x\right)  }\\
&  +\sum_{1\leq m\leq i}\left\|  \frac{\partial^{i}}{\partial x^{i}%
}\triangledown\phi\left(  X,\xi\right)  \right\|  _{L^{\infty}\left(
x\right)  }\left[  \frac{\partial^{j_{m}}}{\partial x^{j_{m}}}X\right]
_{0,\alpha,\left(  x\right)  }\prod_{p\neq m}\left\|  \frac{\partial^{j_{p}}%
}{\partial x^{j_{p}}}X\right\|  _{L^{\infty}\left(  x\right)  }\\
&  \leq\frac{C\left\|  \phi\right\|  _{Y_{\ell,\alpha}}}{\left(  \xi+1\right)
^{N-1}t^{\ell+\alpha}}.
\end{align*}
Putting the above inequalities into (\ref{Gxl}) yields
\begin{align*}
\left\|  \frac{\partial^{\ell}G}{\partial x^{\ell}}\left(  \xi\right)
\right\|  _{L^{\infty}\left(  x\right)  }  &  \leq\frac{C\left\|
\phi\right\|  _{Y_{\ell,\alpha}}}{t^{\ell}\left(  \xi+1\right)  ^{N-1}}%
+\frac{C\left\|  \phi\right\|  _{Y_{\ell,\alpha}}}{\left(  \xi+1\right)  ^{N}%
}\int_{0}^{\xi}\left\|  \frac{\partial^{\ell}\varphi}{\partial x^{\ell}%
}\left(  \bar{s}\right)  \right\|  _{L^{\infty}\left(  x\right)  }d\bar{s}\\
&  \leq\frac{C\left\|  \phi\right\|  _{Y_{\ell,\alpha}}}{t^{\ell}\left(
\xi+1\right)  ^{N-1}}+\frac{C\left\|  \phi\right\|  _{Y_{\ell,\alpha}}%
}{\left(  \xi+1\right)  ^{N-1}}\sup_{0\leq s\leq t}\left\|  \frac
{\partial^{\ell}\varphi}{\partial x^{\ell}}\left(  s\right)  \right\|
_{L^{\infty}\left(  x\right)  },
\end{align*}%
\begin{align*}
\left[  \frac{\partial^{\ell}G}{\partial x^{\ell}}\left(  \xi\right)  \right]
_{0,\alpha,\left(  x\right)  }  &  \leq\frac{C\left\|  \phi\right\|
_{Y_{\ell,\alpha}}}{t^{\ell+\alpha}\left(  \xi+1\right)  ^{N-1}}%
+\frac{C\left\|  \phi\right\|  _{Y_{\ell,\alpha}}}{\left(  \xi+1\right)  ^{N}%
}\int_{0}^{\xi}\left[  \frac{\partial^{\ell}\varphi}{\partial x^{\ell}}\left(
\bar{s}\right)  \right]  _{0,\alpha,\left(  x\right)  }d\bar{s}\\
&  \leq\frac{C\left\|  \phi\right\|  _{Y_{\ell,\alpha}}}{t^{\ell+\alpha
}\left(  \xi+1\right)  ^{N-1}}+\frac{C\left\|  \phi\right\|  _{Y_{\ell,\alpha
}}}{\left(  \xi+1\right)  ^{N-1}}\sup_{0\leq s\leq t}\left[  \frac
{\partial^{\ell}\varphi}{\partial x^{\ell}}\left(  s\right)  \right]
_{0,\alpha,\left(  x\right)  }.
\end{align*}
In a similar manner, we take $\frac{\partial^{\ell+1}}{\partial x^{\ell
}\partial x_{0}}$ of $G$ to get
\begin{align*}
\frac{\partial^{\ell+1}G}{\partial x^{\ell}\partial x_{0}}\left(  \xi\right)
&  =\sum_{\substack{1\leq i\leq\ell,\\j_{1}+...+j_{i}=\ell}}B_{ij_{1}...j_{i}%
}\frac{\partial^{i}}{\partial x^{i}}\triangledown\phi\left(  X,\xi\right)
\frac{\partial\partial^{j_{1}}}{\partial x_{0}\partial x^{j_{1}}}%
X\cdot...\cdot\frac{\partial^{j_{i}}}{\partial x^{j_{i}}}X\\
&  =\triangledown^{2}\phi\left(  X,\xi\right)  \frac{\partial^{\ell+1}%
X}{\partial x^{\ell}\partial x_{0}}+\sum_{\substack{j_{m}<\ell\\1\leq m\leq
i}}...=\triangledown^{2}\phi\left(  X,\xi\right)  \int_{0}^{\xi}\frac
{\partial^{\ell+1}\varphi}{\partial x^{\ell}\partial x_{0}}\left(  \bar
{s}\right)  d\bar{s}+\sum_{_{\substack{j_{m}<\ell\\1\leq m\leq i}}}...,
\end{align*}
We use the induction hypotheses to bound all the terms with $j_{m}<\ell$ for
all $1\leq m\leq i$ on the right-hand side of the above identity as
\begin{align*}
&  \left\|  \frac{\partial^{i}}{\partial x^{i}}\triangledown\phi\left(
X,\xi\right)  \frac{\partial\partial^{j_{1}}}{\partial x_{0}\partial x^{j_{1}%
}}X\cdot...\cdot\frac{\partial^{j_{i}}}{\partial x^{j_{i}}}X\right\|
_{L^{\infty}\left(  x\right)  }\\
&  \leq\left\|  \frac{\partial^{i}}{\partial x^{i}}\triangledown\phi\left(
\xi\right)  \right\|  _{L^{\infty}\left(  x\right)  }\int_{0}^{\xi}\left\|
\frac{\partial\partial^{j_{1}}\varphi}{\partial x_{0}\partial x^{j_{1}}%
}\left(  \bar{s}\right)  \right\|  _{L^{\infty}\left(  x\right)  }d\bar
{s}\cdot...\cdot\left\|  \frac{\partial^{j_{i}}}{\partial x^{j_{i}}}X\right\|
_{L^{\infty}\left(  x\right)  }\\
&  \leq\frac{C\left\|  \phi\right\|  _{Y_{\ell,\alpha}}}{\left(  \xi+1\right)
^{i+N-1}}\frac{C\left\|  \phi\right\|  _{Y_{\ell,\alpha}}}{t^{j_{1}}}%
\frac{C\left\|  \phi\right\|  _{Y_{\ell,\alpha}}\xi}{t^{j_{2}}}\cdot
...\cdot\frac{C\left\|  \phi\right\|  _{Y_{\ell,\alpha}}\xi}{t^{j_{i}}}%
\leq\frac{C\left\|  \phi\right\|  _{Y_{\ell,\alpha}}}{\left(  \xi+1\right)
^{N}t^{\ell}},
\end{align*}%
\begin{align*}
&  \left[  \frac{\partial^{i}}{\partial x^{i}}\triangledown\phi\left(
X,\xi\right)  \frac{\partial\partial^{j_{1}}}{\partial x_{0}\partial x^{j_{1}%
}}X\cdot...\cdot\frac{\partial^{j_{i}}}{\partial x^{j_{i}}}X\right]
_{0,\alpha,\left(  x\right)  }\\
&  \leq\left[  \frac{\partial^{i}}{\partial x^{i}}\triangledown\phi\left(
X,\xi\right)  \right]  _{0,\alpha,\left(  x\right)  }\left(  \int_{0}^{\xi
}\left\|  \frac{\partial\partial^{j_{1}}\varphi}{\partial x_{0}\partial
x^{j_{1}}}\left(  \bar{s}\right)  \right\|  _{L^{\infty}\left(  x\right)
}d\bar{s}\right)  \left(  \prod_{p\neq1}\left\|  \frac{\partial^{j_{p}}%
}{\partial x^{j_{p}}}X\right\|  _{L^{\infty}\left(  x\right)  }\right) \\
&  +\left\|  \frac{\partial^{i}}{\partial x^{i}}\triangledown\phi\left(
X,\xi\right)  \right\|  _{L^{\infty}\left(  x\right)  }\left(  \int_{0}^{\xi
}\left[  \frac{\partial\partial^{j_{1}}\varphi}{\partial x_{0}\partial
x^{j_{1}}}\left(  \bar{s}\right)  \right]  _{0,\alpha,\left(  x\right)  }%
d\bar{s}\right)  \left(  \prod_{p\neq1}\left\|  \frac{\partial^{j_{p}}%
}{\partial x^{j_{p}}}X\right\|  _{L^{\infty}\left(  x\right)  }\right) \\
&  +\sum_{2\leq m\leq i}\left\|  \frac{\partial^{i}}{\partial x^{i}%
}\triangledown\phi\left(  X,\xi\right)  \right\|  _{L^{\infty}\left(
x\right)  }\left(  \int_{0}^{\xi}\left\|  \frac{\partial\partial^{j_{1}%
}\varphi}{\partial x_{0}\partial x^{j_{1}}}\left(  \bar{s}\right)  \right\|
_{L^{\infty}\left(  x\right)  }d\bar{s}\right)  \left(  \left[  \frac
{\partial^{j_{m}}}{\partial x^{j_{m}}}X\right]  _{0,\alpha,\left(  x\right)
}\right)  \left(  \prod_{p\neq m,1}\left\|  \frac{\partial^{j_{p}}}{\partial
x^{j_{p}}}X\right\|  _{L^{\infty}\left(  x\right)  }\right) \\
&  \leq\frac{C\left\|  \phi\right\|  _{Y_{\ell,\alpha}}}{\left(  \xi+1\right)
^{i+N-1+\alpha}}\frac{\xi^{\alpha}}{t^{\alpha}}\frac{C\left\|  \phi\right\|
_{Y_{\ell,\alpha}}}{t^{j_{1}}}\frac{C\left\|  \phi\right\|  _{Y_{\ell,\alpha}%
}\xi}{t^{j_{2}}}\cdot...\cdot\frac{C\left\|  \phi\right\|  _{Y_{\ell,\alpha}%
}\xi}{t^{j_{i}}}\\
&  +\frac{C\left\|  \phi\right\|  _{Y_{\ell,\alpha}}}{\left(  \xi+1\right)
^{i+N-1}}\frac{C\left\|  \phi\right\|  _{Y_{\ell,\alpha}}}{t^{j_{1}+\alpha}%
}\frac{C\left\|  \phi\right\|  _{Y_{\ell,\alpha}}\xi}{t^{j_{2}}}\cdot
...\cdot\frac{C\left\|  \phi\right\|  _{Y_{\ell,\alpha}}\xi}{t^{j_{i}}}\\
&  +\frac{C\left\|  \phi\right\|  _{Y_{\ell,\alpha}}}{\left(  \xi+1\right)
^{i+N-1}}\frac{C\left\|  \phi\right\|  _{Y_{\ell,\alpha}}}{t^{j_{1}}}%
\frac{C\left\|  \phi\right\|  _{Y_{\ell,\alpha}}\xi}{t^{j_{m}+\alpha}}%
\frac{C\left\|  \phi\right\|  _{Y_{\ell,\alpha}}\xi}{t^{j_{2}}}\cdot
...\cdot\frac{C\left\|  \phi\right\|  _{Y_{\ell,\alpha}}\xi}{t^{j_{i}}}\\
&  \leq\frac{C\left\|  \phi\right\|  _{Y_{\ell,\alpha}}}{\left(  \xi+1\right)
^{N}t^{\ell+\alpha}}.
\end{align*}
Thus we obtain
\[
\left\|  \frac{\partial^{\ell}G}{\partial x^{\ell}}\left(  \xi\right)
\right\|  _{L^{\infty}\left(  x\right)  }\leq\frac{C\left\|  \phi\right\|
_{Y_{\ell,\alpha}}}{t^{\ell}\left(  \xi+1\right)  ^{N}}+\frac{C\left\|
\phi\right\|  _{Y_{\ell,\alpha}}}{\left(  \xi+1\right)  ^{N}}\int_{0}^{\xi
}\left\|  \frac{\partial^{\ell+1}\varphi}{\partial x^{\ell}\partial x_{0}%
}\left(  \bar{s}\right)  \right\|  _{L^{\infty}\left(  x\right)  }d\bar{s}%
\]%
\[
\left[  \frac{\partial^{\ell+1}G}{\partial x^{\ell}\partial x_{0}}\left(
\xi\right)  \right]  _{0,\alpha,\left(  x\right)  }\leq\frac{C\left\|
\phi\right\|  _{Y_{\ell,\alpha}}}{t^{\ell+\alpha}\left(  \xi+1\right)  ^{N}%
}+\frac{C\left\|  \phi\right\|  _{Y_{\ell,\alpha}}}{\left(  \xi+1\right)
^{N}}\int_{0}^{\xi}\left[  \frac{\partial^{\ell+1}\varphi}{\partial x^{\ell
}\partial x_{0}}\left(  \bar{s}\right)  \right]  _{0,\alpha,\left(  x\right)
}d\bar{s}.
\]
This completes the proof.
\end{proof}

We now prove Proposition \ref{regulhigh}.

\begin{proof}
[Proof of Proposition \ref{regulhigh}]Taking $\frac{\partial^{\ell}}{\partial
x^{\ell}}$ of $\varphi$ and using (\ref{phi}), we get
\[
\frac{\partial^{\ell}\varphi}{\partial x^{\ell}}\left(  \xi\right)  =-\int
_{s}^{t}\frac{\partial^{\ell}G}{\partial x^{\ell}}\left(  \xi\right)
d\xi+\frac{1}{t}\int_{0}^{t}\xi\frac{\partial^{\ell}G}{\partial x^{\ell}%
}\left(  \xi\right)  d\xi.
\]

$\left|  \left|  \left|  \cdot\right|  \right|  \right|  $ will denote a
generic norm (or seminorm) that will be assumed to take the specific
values $\left\|  \cdot\right\|  _{L^{\infty}},$\ $\left[  \cdot\right]
_{0,\alpha}.\;$\ Then
\begin{equation}
\left|  \left|  \left|  \frac{\partial^{\ell}\varphi}{\partial x^{\ell}%
}\left(  \xi\right)  \right|  \right|  \right|  \leq\int_{s}^{t}\left|
\left|  \left|  \frac{\partial^{\ell}G}{\partial x^{\ell}}\left(  \xi\right)
\right|  \right|  \right|  d\xi+\frac{1}{t}\int_{0}^{t}\xi\left|  \left|
\left|  \frac{\partial^{\ell}G}{\partial x^{\ell}}\left(  \xi\right)  \right|
\right|  \right|  d\xi\leq2\int_{0}^{t}\left|  \left|  \left|  \frac
{\partial^{\ell}G}{\partial x^{\ell}}\left(  \xi\right)  \right|  \right|
\right|  d\xi, \label{phixl}%
\end{equation}
where we used $\frac{\xi}{t}\leq1.$ Suppose that the induction hypothesis
(\ref{hypo1})-(\ref{hypo2}) is satisfied. By Lemma \ref{Lem-Gxl}, we
then\ have
\begin{equation}
\left|  \left|  \left|  \frac{\partial^{\ell}G}{\partial x^{\ell}}\left(
\xi\right)  \right|  \right|  \right|  \leq\frac{C\left\|  \phi\right\|
_{Y_{\ell,\alpha}}}{\left(  \xi+1\right)  ^{N-1}}\{\frac{1}{t^{\gamma}}%
+\sup_{0\leq s\leq t}\left|  \left|  \left|  \frac{\partial^{\ell}\varphi
}{\partial x^{\ell}}\left(  s\right)  \right|  \right|  \right|  \},
\label{Gxl1}%
\end{equation}
where
\begin{equation}
\gamma=\left\{
\begin{array}
[c]{c}%
\ell,\text{ \ \ \ \ \ \ \ \ if }\left|  \left|  \left|  \cdot\right|  \right|
\right|  =\left\|  \cdot\right\|  _{L^{\infty}},\\
\ell+\alpha,\text{ \ \ if }\left|  \left|  \left|  \cdot\right|  \right|
\right|  =[\cdot]_{0,\alpha}%
\end{array}
\right.  . \label{gamma}%
\end{equation}
Putting (\ref{Gxl1}) into (\ref{phixl}) yields
\[
\left|  \left|  \left|  \frac{\partial^{\ell}\varphi}{\partial x^{\ell}%
}\left(  \xi\right)  \right|  \right|  \right|  \leq C\left\|  \phi\right\|
_{Y_{\ell,\alpha}}\left(  \int_{0}^{t}\frac{d\xi}{\left(  \xi+1\right)
^{N-1}}\right)  \{\frac{1}{t^{\gamma}}+\sup_{0\leq s\leq t}\left|  \left|
\left|  \frac{\partial^{\ell}\varphi}{\partial x^{\ell}}\left(  s\right)
\right|  \right|  \right|  \}\leq C\left\|  \phi\right\|  _{Y_{\ell,\alpha}%
}\{\frac{1}{t^{\gamma}}+\sup_{0\leq s\leq t}\left|  \left|  \left|
\frac{\partial^{\ell}\varphi}{\partial x^{\ell}}\left(  s\right)  \right|
\right|  \right|  \}.
\]
Thus if $\varepsilon_{0}$ is small enough, we obtain (\ref{phixl-inf}). In a
similar way, we have
\begin{equation}
\left|  \left|  \left|  \frac{\partial^{\ell+1}\varphi}{\partial x^{\ell
}\partial x_{0}}\left(  s\right)  \right|  \right|  \right|  \leq\int_{s}%
^{t}\left|  \left|  \left|  \frac{\partial^{\ell+1}G}{\partial x^{\ell
}\partial x_{0}}\left(  \xi\right)  \right|  \right|  \right|  d\xi+\frac
{1}{t}\int_{0}^{t}\xi\left|  \left|  \left|  \frac{\partial^{\ell+1}%
G}{\partial x^{\ell}\partial x_{0}}\left(  \xi\right)  \right|  \right|
\right|  d\xi. \label{phixx0l1}%
\end{equation}
Using Lemma \ref{Lem-Gxl} yields
\begin{equation}
\left|  \left|  \left|  \frac{\partial^{\ell+1}G}{\partial x^{\ell}\partial
x_{0}}\left(  \xi\right)  \right|  \right|  \right|  \leq\frac{C\left\|
\phi\right\|  _{Y_{\ell,\alpha}}}{\left(  \xi+1\right)  ^{N}}\{\frac
{1}{t^{\gamma}}+\int_{0}^{\xi}\left|  \left|  \left|  \frac{\partial^{\ell
+1}\varphi}{\partial x^{\ell}\partial x_{0}}\left(  \bar{s}\right)  \right|
\right|  \right|  d\bar{s}\}, \label{Gxx0l}%
\end{equation}
where $\gamma$ is as in (\ref{gamma}). By putting (\ref{Gxx0l}) into
(\ref{phixx0l1}) and integrating from $s=0$ to $s=t$, we have
\begin{align*}
\int_{0}^{t}\left|  \left|  \left|  \frac{\partial^{\ell+1}\varphi}{\partial
x^{\ell}\partial x_{0}}\left(  s\right)  \right|  \right|  \right|  ds  &
\leq C\left\|  \phi\right\|  _{Y_{\ell,\alpha}}\left(  \int_{0}^{t}\int
_{s}^{t}\frac{d\xi}{\left(  \xi+1\right)  ^{N}}ds\right)  \{\frac{1}%
{t^{\gamma}}+\left(  \int_{0}^{t}\left|  \left|  \left|  \frac{\partial
^{\ell+1}\varphi}{\partial x^{\ell}\partial x_{0}}\left(  s\right)  \right|
\right|  \right|  ds\right)  \}\\
&  +C\left\|  \phi\right\|  _{Y_{\ell,\alpha}}\left(  \int_{0}^{t}\frac{\xi
d\xi}{\left(  \xi+1\right)  ^{N}}\right)  \{\frac{1}{t^{\gamma}}+\left(
\int_{0}^{t}\left|  \left|  \left|  \frac{\partial^{\ell+1}\varphi}{\partial
x^{\ell}\partial x_{0}}\left(  s\right)  \right|  \right|  \right|  ds\right)
\}\\
&  \leq C\left\|  \phi\right\|  _{Y_{\ell,\alpha}}\{\frac{1}{t^{\gamma}%
}+\left(  \int_{0}^{t}\left|  \left|  \left|  \frac{\partial^{\ell+1}\varphi
}{\partial x^{\ell}\partial x_{0}}\left(  s\right)  \right|  \right|  \right|
ds\right)  \}.
\end{align*}
Thus we obtain (\ref{phixx0l}) provided $\varepsilon_{0}$ is small. We then
substitute (\ref{Gxx0l}) for (\ref{phixx0l}) and put $s=t$ in (\ref{phixx0l1})
to get
\[
\left|  \left|  \left|  \frac{\partial^{\ell+1}\varphi}{\partial x^{\ell
}\partial x_{0}}\left(  t\right)  \right|  \right|  \right|  \leq\frac{1}%
{t}\frac{C\left\|  \phi\right\|  _{Y_{\ell,\alpha}}}{t^{\gamma}}\int_{0}%
^{t}\frac{\xi d\xi}{\left(  \xi+1\right)  ^{N}}\leq\frac{C\left\|
\phi\right\|  _{Y_{\ell,\alpha}}}{t^{\gamma+1}}.
\]
Therefore the proof is complete.
\end{proof}

\subsection{Estimating the potential $\phi$ in terms of the density $\rho.$}

The following is a standard regularity result for the Poisson equation.

\begin{lemma}
\label{elliptic}(Elliptic regularity theory)
\[
\left\|  \phi\right\|  _{Y_{k,\alpha}}\leq C\left\|  \rho\right\|
_{X_{k,\alpha}}%
\]

\end{lemma}

\begin{proof}
For any fixed $t>0$ we define
\[
\tilde{\rho}\left(  z,t\right)  =\left(  t+1\right)  ^{N}\rho\left(  z\left(
t+1\right)  ,t\right)
\]

Notice that:
\[
\int\left|  \tilde{\rho}\left(  z\right)  \right|  d^{N}z=\int\left(
t+1\right)  ^{N}\left|  \rho\left(  z\left(  t+1\right)  ,t\right)  \right|
d^{N}z=\int\left|  \rho\left(  x\right)  \right|  d^{N}x
\]

On the other hand
\begin{align*}
\sum_{\ell=0}^{k}\left\Vert \nabla_{z}^{\ell}\tilde{\rho}\left(
\cdot,t\right)  \right\Vert _{L^{\infty}\left(  \mathbb{R}^{N}\right)  }  &
=\left(  t+1\right)  ^{N}\sum_{\ell=0}^{k}\left(  t+1\right)  ^{\ell
}\left\Vert \nabla_{x}^{\ell}\rho\left(  \cdot,t\right)  \right\Vert
_{L^{\infty}\left(  \mathbb{R}^{N}\right)  }\\
\sup_{z_{1},z_{2}\in\mathbb{R}^{N}}\frac{\left\vert \nabla_{z}^{k}\tilde{\rho
}\left(  z_{1},t\right)  -\nabla_{z}^{k}\tilde{\rho}\left(  z_{2},t\right)
\right\vert }{\left\vert z_{1}-z_{2}\right\vert ^{\alpha}}  &  =\left(
t+1\right)  ^{N+k+\alpha}\sup_{x_{1},x_{2}\in\mathbb{R}^{N}}\frac{\left\vert
\nabla_{x}^{k}\rho\left(  x_{1},t\right)  -\nabla_{x}^{k}\rho\left(
x_{2},t\right)  \right\vert }{\left\vert x_{1}-x_{2}\right\vert ^{\alpha}}%
\end{align*}

Then:
\begin{equation}
\int\left|  \tilde{\rho}\left(  z\right)  \right|  d^{3}z+\sum_{\ell=0}%
^{k}\left\|  \nabla_{z}^{\ell}\tilde{\rho}\left(  \cdot,t\right)  \right\|
_{L^{\infty}\left(  \mathbb{R}^{N}\right)  }+\sup_{z_{1},z_{2}\in
\mathbb{R}^{N}}\frac{\left|  \nabla_{z}^{k}\tilde{\rho}\left(  z_{1},t\right)
-\nabla_{z}^{k}\tilde{\rho}\left(  z_{2},t\right)  \right|  }{\left|
z_{1}-z_{2}\right|  ^{\alpha}}\leq C\left\|  \rho\right\|  _{X_{k,\alpha},}
\label{E1}%
\end{equation}

On the other hand, by assumption
\[
\Delta_{x}\phi=\rho
\]

We define
\[
\tilde{\phi}\left(  z\right)  =\left(  t+1\right)  ^{N-2}\phi\left(  z\left(
t+1\right)  ,t\right)
\]

Then:
\[
\Delta_{z}\tilde{\phi}=\tilde{\rho}%
\]

We now claim that the following estimate holds
\begin{align}
&  \sum_{\ell=1}^{k+2}\left\Vert \nabla_{z}^{\ell}\tilde{\phi}\left(
\cdot,t\right)  \right\Vert _{L^{\infty}\left(  \mathbb{R}^{N}\right)  }%
+\sup_{z_{1},z_{2}\in\mathbb{R}^{N}}\frac{\left\vert \nabla_{z}^{k+2}%
\tilde{\phi}\left(  z_{1},t\right)  -\nabla_{z}^{k+2}\tilde{\phi}\left(
z_{2},t\right)  \right\vert }{\left\vert z_{1}-z_{2}\right\vert ^{\alpha}}\leq
CJ\label{Poissreg}\\
J  &  \equiv\left[  \int_{\mathbb{R}^{N}}\left\vert \tilde{\rho}\left(
z\right)  \right\vert d^{N}z+\sum_{\ell=0}^{k}\left\Vert \nabla_{z}^{\ell
}\tilde{\rho}\left(  \cdot,t\right)  \right\Vert _{L^{\infty}\left(
\mathbb{R}^{N}\right)  }+\sup_{z_{1},z_{2}\in\mathbb{R}^{N}}\frac{\left\vert
\nabla_{z}^{k}\tilde{\rho}\left(  z_{1},t\right)  -\nabla_{z}^{k}\tilde{\rho
}\left(  z_{2},t\right)  \right\vert }{\left\vert z_{1}-z_{2}\right\vert
^{\alpha}}\right] \nonumber
\end{align}

Indeed, a standard interpolation argument yields%
\[
\left\Vert \tilde{\rho}\right\Vert _{L^{p}}\leq\left\Vert \tilde{\rho
}\right\Vert _{L^{\infty}}^{\frac{p-1}{p}}\left\Vert \tilde{\rho}\right\Vert
_{L^{1}}^{\frac{1}{p}}\leq J\;\;,\;\;1\leq p\leq\infty
\]

Using then the Calderon-Zygmund inequality it follows that:%
\[
\left\|  \nabla_{z}^{2}\tilde{\phi}\right\|  _{L^{p}}\leq C\left\|
\tilde{\rho}\right\|  _{L^{p}}\leq CJ\;,\;\;1<p<\infty
\]
Therefore, the Sobolev embedding theorem implies that%
\[
\left\Vert \tilde{\phi}\right\Vert_{L^{q}}
\leq CJ\;,\;\;\frac{N}{N-2}<q<\infty
\]
Interior estimates for the Poisson equation in Sobolev spaces give a uniform
bound on the $W^{k+2,q}$ norm of the restriction of $\tilde\phi$ to any unit
ball and hence of the $C^\alpha$ norm of this restriction. 
Using this estimate, (\ref{Poissreg}) follows from the inequality%
\[
\left\|  \tilde{\phi}\right\|  _{k+2,\alpha;B_{\frac{1}{2}}\left(
x_{0}\right)  }\leq C\left[  \left\|  \tilde{\phi}\right\|  _{C^{\alpha}\left(
B_{1}\left(  x_{0}\right)  \right)  }+\left\|  \tilde{\rho}\right\|
_{k,\alpha;B_{1}\left(  x_{0}\right)  }\right]  \leq CJ
\]
that is just a consequence of classical interior estimates for the Poisson
equation (cf. \cite{GT}).
Using the estimate (\ref{E1}) it can be concluded that
\[
\sum_{\ell=1}^{k+2}\left\Vert \nabla_{z}^{\ell}\tilde{\phi}\left(
\cdot,t\right)  \right\Vert _{L^{\infty}\left(  \mathbb{R}^{N}\right)  }%
+\sup_{z_{1},z_{2}\in\mathbb{R}^{N}}\frac{\left\vert \nabla_{z}^{k+2}%
\tilde{\phi}\left(  z_{1},t\right)  -\nabla_{z}^{k+2}\tilde{\phi}\left(
z_{2},t\right)  \right\vert }{\left\vert z_{1}-z_{2}\right\vert ^{\alpha}}\leq
C\left\Vert \rho\right\Vert _{X_{k,\alpha}}%
\]

Using the definition of $\tilde{\phi}$ as well as the definition of the norm
$\left\Vert \cdot\right\Vert _{Y_{k,\alpha}}$ as in (\ref{Yka})$,$ it then
follows that
\[
\left\Vert \phi\right\Vert _{Y_{k,\alpha}}\leq C\left\Vert \rho\right\Vert
_{X_{k,\alpha}}%
\]
and this completes the proof of the lemma.
\end{proof}

\subsection{Conservation of the $L^{1}$ norm of $f$.}

The following result is standard in the theory of the Vlasov-Poisson equation.
See for instance \cite{glassey1}.

\begin{lemma}
Suppose that $f\left(  x,v,t\right)  $ solves the problem (\ref{VP1}),
(\ref{VP2}). Then:
\begin{equation}
\int_{\mathbb{R}^{N}}\int_{\mathbb{R}^{N}}\left|  f\left(  x,v,t\right)
\right|  dvdx=\int_{\mathbb{R}^{N}}\int_{\mathbb{R}^{N}}\left|  f_{0}\left(
x,v\right)  \right|  dvdx\;\;,\;\;t>0 \label{masscon}%
\end{equation}

\end{lemma}

\subsection{The proof of Theorem \ref{thm}.\bigskip}

We now prove our main Theorem \ref{thm}. We use a 
continuation argument as well as the assumptions on $f_{0}$
to derive estimates for $\rho$ and to 
close the argument. More precisely, we will show that an
estimate of the form%
\begin{align}
&  \int\left\vert \rho\left(  x,t\right)  \right\vert d^{N}x+\left(
t+1\right)  ^{N}\sum_{\ell=0}^{k}\left(  t+1\right)  ^{\ell}\left\Vert
\nabla^{\ell}\rho\left(  \cdot,t\right)  \right\Vert _{L^{\infty}\left(
\mathbb{R}^{N}\right)  }\nonumber\\
&  +\left(  t+1\right)  ^{N+k+\alpha}\sup_{x,y\in\mathbb{R}^{N}}%
\frac{\left\vert \nabla^{k}\rho\left(  x,t\right)  -\nabla^{k}\rho\left(
y,t\right)  \right\vert }{\left\vert x-y\right\vert ^{\alpha}}\leq
M\varepsilon_{0}\nonumber
\end{align}
for $0\leq t\leq t^{\ast}$ implies an estimate of the form%
\begin{align}
&  \int\left\vert \rho\left(  x,t\right)  \right\vert d^{N}x+\left(
t+1\right)  ^{N}\sum_{\ell=0}^{k}\left(  t+1\right)  ^{\ell}\left\Vert
\nabla^{\ell}\rho\left(  \cdot,t\right)  \right\Vert _{L^{\infty}\left(
\mathbb{R}^{N}\right)  }\nonumber\\
&  +\left(  t+1\right)  ^{N+k+\alpha}\sup_{x,y\in\mathbb{R}^{N}}%
\frac{\left\vert \nabla^{k}\rho\left(  x,t\right)  -\nabla^{k}\rho\left(
y,t\right)  \right\vert }{\left\vert x-y\right\vert ^{\alpha}}\leq
C\varepsilon_{0}\label{bound2}%
\end{align}
for $0\leq t\leq t^{\ast}+\delta\left(  t^{\ast}\right)  ,$ where
$\delta\left(  t^{\ast}\right)  >0$ and $C$ is independent of $M$ if
$\varepsilon_{0}$ is small enough. 
Together with a suitable local existence theorem which guarantees that an 
inequality of the form (\ref{bound2}) holds on some time interval $[0,t_1]$ 
for $t_1>0$ the estimates which have been derived imply that (\ref{bound2}) 
can be extended to the whole time interval $0\le t<\infty$ and this implies 
Theorem 1. Notice that this in particular yields a global existence theorem 
generalizing that of \cite{bardos} for $N=3$. The local existence theorem can 
be obtained by combining the estimates obtained in this paper with a 
contraction mapping argument in a straightforward way.

Notice that the decay assumptions on $f_{0}$ have not been used until now and
only the decay properties of the potential $\phi$ have been used. We now use
the decay properties of $f_{0}$ for the first time in the following proof.

\begin{proof}
Suppose first that $0\leq t\leq1.$ Let $\left(  x_{0},v_{0}\right)  $ denote
the starting point for the solution of the characteristic equations reaching
the point $\left(  x,v\right)  $ at time $t.$ More precisely,
\begin{equation}
\frac{d\bar{X}\left(  s\right)  }{ds}=\bar{V}\left(  s\right)  ,\ \frac
{d\bar{V}\left(  s\right)  }{ds}=\nabla\phi\left(  \bar{X}\left(  s\right)
,s\right)  ,\ \bar{X}\left(  t\right)  =x,\ \bar{V}\left(  t\right)  =v
\label{Xchar}%
\end{equation}
Notice that $\bar{X}\left(  s\right)  =\bar{X}\left(  s;x,v,t\right)
,\;\bar{V}\left(  s\right)  =\bar{V}\left(  s;x,v,t\right)  $ and:
\[
x_{0}=x_{0}\left(  x,v,t\right)  =\bar{X}\left(  0;x,v,t\right)
,\;\;\;v_{0}=v_{0}\left(  x,v,t\right)  =\bar{V}\left(  0;x,v,t\right)  .
\]
Then
\begin{equation}
\left|  v_{0}\right|  =\left|  v_{0}\left(  x,v,t\right)  \right|  \geq\left|
v\right|  -\int_{0}^{t}\left|  \nabla\phi\left(  s\right)  \right|
ds\geq\left|  v\right|  -C\left\|  \phi\right\|  _{Y_{k,\alpha}}\geq\left|
v\right|  -\frac{1}{2}, \label{v}%
\end{equation}
if $\varepsilon_{0}$ is small enough. Taking the derivative $\frac
{\partial^{\ell}}{\partial x^{\ell}}$ of the formula
\[
\rho\left(  x,t\right)  =\int f_{0}\left(  x_{0},v_{0}\right)  dv
\]
yields, for $0\leq\ell\leq k,$%
\[
\frac{\partial^{\ell}\rho}{\partial x^{\ell}}\left(  x,t\right)
=\sum_{\substack{j_{1}+...+j_{i}=\ell,\\0\leq i\leq\ell,\ 0\leq m\leq
i}}C_{ijpm}\int\frac{\partial^{i}f_{0}}{\partial x^{i-m}\partial v^{m}}\left(
x_{0},v_{0}\right)  \frac{\partial^{j_{1}}v_{0}}{\partial x^{j_{1}}}%
\cdot...\cdot\frac{\partial^{j_{i}}x_{0}}{\partial x^{j_{i}}}dv,
\]
Notice that the derivatives $\frac{\partial^{j_{1}}v_{0}}{\partial x^{j_{1}}%
},...,\frac{\partial^{j_{i}}x_{0}}{\partial x^{j_{i}}}$\ are bounded $0\leq
t\leq1$\ as $C\left(  1+\left\|  \phi\right\|  _{Y_{\ell,\alpha}}\right)
,$\ as can be seen by differentiating the characteristic equations
(\ref{Xchar}) with respect to $x$\ and $v.$ It is then straightforward to see
that for $0\leq\ell\leq k,$ using (\ref{smalldata}) and (\ref{v}) yields
\begin{align}
\left\|  \frac{\partial^{\ell}\rho}{\partial x^{\ell}}\left(  t\right)
\right\|  _{L^{\infty}\left(  x\right)  }  &  \leq C\left(  1+\left\|
\phi\right\|  _{Y_{\ell,\alpha}}\right)  \int\left|  \frac{\partial^{i}f_{0}%
}{\partial x^{i-m}\partial v^{m}}\left(  x_{0},v_{0}\right)  \right|
dv\label{rhoder1}\\
&  \leq C\varepsilon_{0}\left(  1+\left\|  \phi\right\|  _{Y_{\ell,\alpha}%
}\right)  \int\frac{dv}{\left(  1+\left|  v\right|  \right)  ^{K}}\nonumber\\
&  \leq C_{\ell}\varepsilon_{0}\left(  1+\left\|  \phi\right\|  _{Y_{\ell
,\alpha}}\right)  \;,\;\;0\leq t\leq1\nonumber
\end{align}%
\begin{align}
\left[  \frac{\partial^{\ell}\rho}{\partial x^{\ell}}\left(  t\right)
\right]  _{0,\alpha,\left(  x\right)  }  &  \leq C\varepsilon_{0}\left(
1+\left\|  \phi\right\|  _{Y_{\ell,\alpha}}\right)  \int\frac{dv}{\left(
1+\left|  v\right|  \right)  ^{K}}\label{rhoder2}\\
&  \leq C_{\ell}\varepsilon_{0}\left(  1+\left\|  \phi\right\|  _{Y_{\ell
,\alpha}}\right)  \;,\;\;0\leq t\leq1\nonumber
\end{align}
Next we treat the case $t\geq1.$ By taking $\frac{\partial^{\ell}%
}{\partial x^{\ell}}$ of $\rho\left(  x,t\right)  $ in (\ref{rho}), we get,
for $0\leq\ell\leq k,$%
\[
\frac{\partial^{\ell}\rho}{\partial x^{\ell}}\left(  x,t\right)
=\sum_{\substack{j_{1}+...+j_{i}+p=\ell\\0\leq i\leq\ell}}C_{ijp}\int
\frac{\partial^{i}f_{0}}{\partial v^{i}}\left(  x_{0},V\left(  0\right)
\right)  \left(  \frac{\partial}{\partial x}\right)  ^{j_{1}}V\left(
0\right)  \cdot...\cdot\left(  \frac{\partial}{\partial x}\right)  ^{j_{i}%
}V\left(  0\right)  \frac{\partial^{p}}{\partial x^{p}}\left(  \left|
\det\frac{\partial w}{\partial x_{0}}\right|  \right)  dx_{0},
\]
where
\[
V\left(  0\right)  =V\left(  0,t,x,w\left(  t,x,x_{0}\right)  \right)
=\frac{x-x_{0}}{t}+\varphi\left(  0;t,x,x_{0}\right)  .
\]
By Lemma \ref{phidecayk0}, Lemma \ref{phidecayk1} and the assumption
(\ref{smalldata}), it is easy to see that for $0\leq\ell\leq k,$%
\begin{align}
\left|  \frac{\partial^{\ell}\rho}{\partial x^{\ell}}\left(  x,t\right)
\right|   &  \leq C\sum_{\substack{j_{1}+...+j_{i}+p=\ell\\0\leq i\leq\ell
}}\frac{\left(  1+C\left\|  \phi\right\|  _{Y_{\ell,\alpha}}\right)
}{t^{j_{1}+...+j_{i}+p+N}}\int\left|  \frac{\partial^{i}f_{0}}{\partial v^{i}%
}\left(  x_{0},V\left(  0\right)  \right)  \right|  dx_{0}\label{rhoder3}\\
&  \leq\frac{C\varepsilon_{0}\left(  1+\left\|  \phi\right\|  _{Y_{\ell
,\alpha}}\right)  }{t^{\ell+N}}\int\frac{dx_{0}}{\left(  1+\left|
x_{0}\right|  \right)  ^{K}}\leq\frac{C\varepsilon_{0}\left(  1+\left\|
\phi\right\|  _{Y_{\ell,\alpha}}\right)  }{t^{\ell+N}}\;,\;\;t>1\nonumber
\end{align}
Using Lemma \ref{Schauder} with $\ell=k$ yields
\begin{align}
\left[  \frac{\partial^{k}\rho}{\partial x^{k}}\left(  t\right)  \right]
_{0,\alpha,\left(  x\right)  }  &  \leq C\frac{\left(  1+C\left\|
\phi\right\|  _{Y_{k,\alpha}}\right)  }{t^{j_{1}+...+j_{i}+p+N+\alpha}}%
\int\left\{  \left|  \frac{\partial^{i}f_{0}}{\partial v^{i}}\left(
x_{0},V\left(  0\right)  \right)  \right|  +\left[  \frac{\partial^{i}f_{0}%
}{\partial v^{i}}\left(  x_{0},\cdot\right)  \right]  _{0,\alpha,\left(
v\right)  }\right\}  dx_{0}\label{rhoder4}\\
&  \leq\frac{C\varepsilon_{0}\left(  1+\left\|  \phi\right\|  _{Y_{k,\alpha}%
}\right)  }{t^{k+N+\alpha}}\int\frac{dx_{0}}{\left(  1+\left|  x_{0}\right|
\right)  ^{K}}\leq\frac{C\varepsilon_{0}\left(  1+\left\|  \phi\right\|
_{Y_{k,\alpha}}\right)  }{t^{k+N+\alpha}}\;,\;\;t>1\nonumber
\end{align}

On the other hand, combining (\ref{masscon}) with the decay assumptions for
$f_{0}$ in (\ref{smalldata}), it follows that
\begin{equation}
\int\left|  \rho\left(  x,t\right)  \right|  dx\leq C\varepsilon_{0}
\label{rhomass}%
\end{equation}

It then follows from (\ref{rhoder1})-(\ref{rhomass}), as well as from Lemma
\ref{elliptic} that
\[
\left\|  \rho\right\|  _{X_{k,\alpha}}\leq C\varepsilon_{0}+C\varepsilon
_{0}\left\|  \phi\right\|  _{Y_{k,\alpha}}\leq C\varepsilon_{0}+C\varepsilon
_{0}\left\|  \rho\right\|  _{X_{k,\alpha}}.
\]
Choosing $\varepsilon_{0}$ small enough, Theorem \ref{thm} follows.
\end{proof}

\bigskip

\section{CONVERGENCE TO THE SELF-SIMILAR BEHAVIOUR.}

We define the following set of\textbf{\ }self-similar variables
\begin{equation}
f\left(  x,v,t\right)  =\frac{1}{\left(  t+1\right)  ^{N}}g\left(
y,v,\tau\right)  ,\label{SS1}%
\end{equation}%
\begin{equation}
\phi\left(  x,v,t\right)  =\frac{1}{(t+1)^{N-2}}\Phi\left(  y,v,\tau\right)
,\label{SS2}%
\end{equation}
where
\begin{equation}
\ \ y=\frac{x}{\left(  t+1\right)  }\;\;,\;\;\tau=\log\left(  t+1\right)
\;.\;\label{SS3}%
\end{equation}
A straightforward computation yields the following transformed system
\begin{equation}
g_{\tau}+\left(  v-y\right)  \cdot\nabla_{y}g+e^{-\left(  N-2\right)  \tau
}\nabla_{y}\Phi\cdot\nabla_{v}g=Ng,\label{g}%
\end{equation}%
\begin{equation}
\Delta_{y}\Phi=\int g\left(  y,v,\tau\right)  dv\equiv\bar{\rho}\left(
y,\tau\right)  ,\label{Phirho}%
\end{equation}
where $g\left(  x,v,0\right)  =g_{0}\left(  x,v\right)  =f_{0}\left(
x,v\right)  =f\left(  x,v,0\right)  .$%

\[
\left\Vert \bar{\rho}\right\Vert _{X_{k,\alpha}}=\sup_{t\geq0}\left\{
\int_{\mathbb{R}^{N}}\left\vert \bar{\rho}\left(  y,t\right)  \right\vert
dy+\sum_{\ell=0}^{k}\left\Vert \nabla^{\ell}\bar{\rho}\left(  \cdot,t\right)
\right\Vert _{L^{\infty}\left(  \mathbb{R}^{N}\right)  }+\sup_{y,y^{\prime}%
\in\mathbb{R}^{N}}\frac{\left\vert \nabla^{k}\bar{\rho}\left(  y,t\right)
-\nabla^{k}\bar{\rho}\left(  y^{\prime},t\right)  \right\vert }{\left\vert
y-y^{\prime}\right\vert ^{\alpha}},\ 0<\alpha<1\right\}  ,
\]

\[
\left\Vert \Phi\right\Vert _{Y_{k,\alpha}}=\sup_{t\geq0}\left\{  \sum_{\ell
=1}^{k+2}\left\Vert \nabla^{\ell}\Phi\left(  \cdot,t\right)  \right\Vert
_{L^{\infty}\left(  \mathbb{R}^{N}\right)  }+\sup_{y,y^{\prime}\in
\mathbb{R}^{N}}\frac{\left\vert \nabla^{k+2}\Phi\left(  y,t\right)
-\nabla^{k+2}\Phi\left(  y^{\prime},t\right)  \right\vert }{\left\vert
y-y^{\prime}\right\vert ^{\alpha}}\right\}  .
\]

\bigskip Notice that Lemma \ref{elliptic} is also valid in self-similar variables

\begin{lemma}
\label{ellipticSS}(Elliptic regularity theory)
\[
\left\|  \Phi\right\|  _{Y_{k,\alpha}}\leq C\left\|  \bar{\rho}\right\|
_{X_{k,\alpha}}%
\]

\end{lemma}

We reformulate Theorem \ref{thm} in self-similar variables

\begin{theorem}
\label{thmSS}Suppose that $g_{0}\left(  y,v\right)  $ satisfies the following
estimates
\begin{align}
\sum_{\ell=0}^{k}\sum_{m=0}^{\ell}\left\vert \frac{\partial^{\ell}g_{0}%
}{\partial x^{m}\partial v^{\ell-m}}\right\vert  &  \leq\frac{\varepsilon_{0}%
}{\left(  1+\left\vert x\right\vert \right)  ^{K}\left(  1+\left\vert
v\right\vert \right)  ^{K}}\;\;,\\
\sum_{m=0}^{k}\sup_{x,x^{\prime}\in IR^{N}}\frac{\left\vert \frac{\partial
^{k}g_{0}}{\partial x^{m}\partial v^{k-m}}\left(  x,v\right)  -\frac
{\partial^{k}g_{0}}{\partial x^{m}\partial v^{k-m}}\left(  x^{\prime
},v\right)  \right\vert }{\left\vert x-x^{\prime}\right\vert ^{\alpha}} &
\leq\frac{\varepsilon_{0}}{\left(  1+\left\vert v\right\vert \right)  ^{K}%
}\;,\;0<\alpha<1\;,\nonumber\\
\sum_{m=0}^{k}\sup_{\left\vert v^{\prime}-v\right\vert \leq1}\frac{\left\vert
\frac{\partial^{k}g_{0}}{\partial x^{m}\partial v^{k-m}}\left(  x,v\right)
-\frac{\partial^{k}g_{0}}{\partial x^{m}\partial v^{k-m}}\left(  x,v^{\prime
}\right)  \right\vert }{\left\vert v-v^{\prime}\right\vert ^{\alpha}} &
\leq\frac{\varepsilon_{0}}{\left(  1+\left\vert x\right\vert \right)
^{K}\left(  1+\left\vert v\right\vert \right)  ^{K}}\;,\;0<\alpha<1\;,
\end{align}
where $K>N$ and $\varepsilon_{0}$ is small enough. Then there exists a 
corresponding solution of the rescaled Vlasov-Poisson system with
\[
\left\Vert \bar{\rho}\right\Vert _{X_{k,\alpha}}\leq C_{k}\varepsilon_{0}%
\]

\end{theorem}

\bigskip The main theorem that we prove in this Section is the following

\begin{theorem}
\label{asymp}Suppose that the assumptions of Theorem \ref{thmSS} are
satisfied. Then, there exist $g_{\infty}\left(  y,y_{0}\right)  \in
C_{loc}^{k},\ \bar{\rho}_{\infty}\left(  y\right)  \in C_{loc}^{k}\cap
L^{1}\left(  \mathbb{R}^{N}\right)  ,\ \Phi_{\infty}\left(  y\right)  \in
C_{loc}^{k+1,\beta}$ satisfying
\begin{align}
e^{-N\tau}g\left(  y,y_{0},\tau\right)   &  \rightarrow g_{\infty}\left(
y,y_{0}\right)  ,\ \ \text{in }C_{loc}^{k}\label{asymp g}\\
\bar{\rho}\left(  y,\tau\right)   &  \rightarrow\bar{\rho}_{\infty}\left(
y\right)  ,\ \ \ \ \ \ \text{in }C_{loc}^{k}\nonumber\\
\Phi\left(  y,\tau\right)   &  \rightarrow\Phi_{\infty}\left(  y\right)
,\ \ \ \ \ \text{in }C_{loc}^{k+1,\beta},\nonumber
\end{align}
for any $0<\beta<1,$ as $\tau\rightarrow\infty.$ Moreover, we have%
\[
\left\Vert \bar{\rho}_{\infty}\right\Vert _{L^{1}\left(  \mathbb{R}%
^{N}\right)  }=\left\Vert g_{0}\right\Vert _{L^{1}\left(  \mathbb{R}^{N}%
\times\mathbb{R}^{N}\right)  }%
\]
and we have the following representation formulae for $g_{\infty}$
\begin{align}
g_{\infty}\left(  y,y_{0}\right)   &  =g_{0}\left(  y_{0},y+\omega_{\infty
}\left(  0,y,y_{0}\right)  \right)  \nonumber\\
\bar{\rho}_{\infty}\left(  y\right)   &  =\int g_{\infty}\left(
y,y_{0}\right)  J_{\infty}\left(  y,y_{0}\right)  dy_{0}\label{S4E2}\\
\Delta_{y}\Phi_{\infty}\left(  y\right)   &  =\bar{\rho}_{\infty}\left(
y\right)  \nonumber
\end{align}
as well as the limit formula, as $\tau\rightarrow\infty$%
\begin{equation}
g\left(  y,v,\tau\right)  \rightarrow\int g_{\infty}\left(  y_{0},y\right)
\delta\left(  v-y\right)  J_{\infty}\left(  y,y_{0}\right)  dy_{0}%
\ ,\;\;\;\text{in\ }\mathcal{D}^{\prime}\left(  \mathbb{R}^{N}\times
\mathbb{R}^{N}\right)  \label{S4E1}%
\end{equation}
where $\omega_{\infty}\left(  s;y,y_{0}\right)  $ is the solution of the
following integral equation
\begin{equation}
\omega_{\infty}\left(  s;y,y_{0}\right)  =-\int_{s}^{\infty}e^{-\left(
N-2\right)  \xi}\nabla_{y}\Phi\left(  y+\left(  y_{0}-y\right)  e^{-\xi}%
+\int_{0}^{\xi}e^{-\left(  \xi-\eta\right)  }\omega_{\infty}\left(
\eta;y,y_{0}\right)  d\eta,\xi\right)  d\xi,\label{omega_inf}%
\end{equation}
and where $J_{\infty}\left(  y,y_{0}\right)  $ is given by
\[
J_{\infty}\left(  y,y_{0}\right)  =\lim_{\tau\rightarrow\infty}\left\vert
\det\left(  -I_{N}+e^{\tau}\frac{\partial\omega_{\infty}}{\partial y_{0}%
}\left(  \tau;y,y_{0}\right)  \right)  \right\vert
\]

\end{theorem}

\begin{remark}
\bigskip Notice that (\ref{asymp g}) can be read in the original set of
variables as
\[
\rho\left(  x,t\right)  \sim\frac{1}{t^{N}}\bar{\rho}_{\infty}\left(  \frac
{x}{\left(  t+1\right)  }\right)  +o\left(  \frac{1}{t^{N}}\right)
\]
as $t\rightarrow\infty,$ uniformly on sets $\left\vert x\right\vert \leq Ct.$
\end{remark}

\begin{remark}
The function $\omega_{\infty}\left(  s;y,y_{0}\right)  $ is small for small
densities. In particular the representation formula (\ref{S4E2}) implies that
the rescaled density function $\bar{\rho}_{\infty}\left(  y\right)  $
approaches the one associated to the free streaming case, defined in
(\ref{S0E1}) if $\varepsilon_{0}\rightarrow0.$ Notice that this shows that the
particular profile that describes the self-similar behaviour of the solutions
depends very sensitively on the initial data $g_{0}.$ This contrasts with the
situation in the one-dimensional case where the leading self-similar behaviour
depends only on the mass of the initial distribution but it does not depend on
any other information on the initial data $g_{0}$ (cf. \cite{bkr}). However,
notice that it is not possible to obtain a closed form expression for 
$g_{\infty}$ in terms of $g_{0}$ due to the fact that the function 
$\omega_{\infty}\left(s;y,y_{0}\right)  $ depends on the values of the 
function $\Phi$ for any $t\in\left(  0,\infty\right)  .$
\end{remark}

In order to prove Theorem \ref{asymp}, we introduce some changes of variables
analogous to the ones used in the previous Section.

Suppose that the characteristics starting at $y_{0},$ $v_{0}$ reach the points
$y,$ $v$ at time $\tau$ and we regard $v_{0}=w_{0}\left(  y,y_{0},\tau\right)
,$ $v=w\left(  y,y_{0},\tau\right)  $ as functions of $y,y_{0},$ and $\tau$
and make the change of variables from $v$ to $y_{0}$ to get
\[
dv=\left|  \det\frac{\partial w}{\partial y_{0}}\right|  dy_{0}%
\]%
\[
\rho\left(  y,\tau\right)  =\int_{\mathbb{R}^{3}}g\left(  y,v,\tau\right)
dv=\int e^{N\tau}g_{0}\left(  y_{0},v_{0}\right)  dv=\int e^{N\tau}%
g_{0}\left(  y_{0},w_{0}\left(  y,y_{0},t\right)  \right)  \left|  \det
\frac{\partial v}{\partial y_{0}}\right|  dy_{0}.
\]
The corresponding boundary value problem in the self similar variables
$\left(  y,v,\tau\right)  $ reads
\begin{align*}
\frac{dY}{ds}  &  =V-Y,\ \frac{dV}{ds}=e^{-\left(  N-2\right)  s}\nabla
_{y}\Phi\left(  Y\left(  s\right)  ,s\right)  ,\ \frac{dg}{ds}=Ng\\
Y\left(  \tau\right)   &  =y,\ Y\left(  0\right)  =y_{0}.
\end{align*}
In the absence of the field, we solve
\begin{align*}
\frac{d\tilde{Y}}{ds}  &  =\tilde{V}-\tilde{Y},\ \frac{d\tilde{V}}{ds}=0,\ \\
\tilde{Y}\left(  \tau\right)   &  =y,\ \tilde{Y}\left(  0\right)  =y_{0},
\end{align*}
which yields
\[
\tilde{V}\left(  s\right)  =\frac{y-y_{0}e^{-\tau}}{1-e^{-\tau}},\ \tilde
{Y}\left(  s\right)  =\frac{y-y_{0}e^{-\tau}}{1-e^{-\tau}}+\frac{y_{0}%
-y}{1-e^{-\tau}}e^{-s}.
\]
As in the previous section, we formulate the above as a perturbed problem from
the free streaming one.
\[
V\equiv\tilde{V}+\omega,\ Y\equiv\tilde{Y}+\zeta,
\]%
\begin{align*}
\frac{d\zeta}{ds}  &  =\omega-\zeta,\ \frac{d\omega}{ds}=e^{-\left(
N-2\right)  s}\nabla_{y}\Phi\left(  \tilde{Y}+\zeta,s\right)  ,\\
\zeta\left(  \tau\right)   &  =\zeta\left(  0\right)  =0.
\end{align*}
It is straightforward to see that
\begin{equation}
\omega\left(  s\right)  =-\int_{s}^{\tau}e^{-\left(  N-2\right)  \xi}%
\nabla_{y}\Phi\left(  Y\left(  \xi\right)  ,\xi\right)  d\xi+\frac{e^{-\tau}%
}{1-e^{-\tau}}\int_{0}^{\tau}e^{-\left(  N-3\right)  \xi}\left(  1-e^{-\xi
}\right)  \nabla_{y}\Phi\left(  Y\left(  \xi\right)  ,\xi\right)  d\xi,
\label{omega}%
\end{equation}%
\[
\zeta\left(  s\right)  =\int_{0}^{s}e^{-\left(  s-\xi\right)  }\omega\left(
\xi\right)  d\xi,
\]
where
\[
Y\left(  \xi\right)  =\frac{y-y_{0}e^{-\tau}}{1-e^{-\tau}}+\frac{y_{0}%
-y}{1-e^{-\tau}}e^{-\xi}+\int_{0}^{\xi}e^{-\left(  \xi-\eta\right)  }%
\omega\left(  \eta\right)  d\eta.
\]
Along the characteristics, we have
\[
\frac{\partial v}{\partial y_{0}}=-\left(  \frac{1}{1-e^{-\tau}}\right)
e^{-\tau}I_{N}+\frac{\partial\omega}{\partial y_{0}}\left(  t\right)
\]
The following result provides some decay estimates for the derivatives of
$\omega,$ analogous to the ones derived in Lemma \ref{phidecayk0}.

\begin{lemma}
\label{omega0}There exists $\varepsilon_{0}$ small such that for any $\tau
\geq1$ and any function $\Phi$ satisfying
\[
\left\|  \Phi\right\|  _{Y_{0,\alpha}}\leq\varepsilon_{0},
\]
we have
\begin{equation}
\int_{0}^{\tau}\left\|  \frac{\partial\omega}{\partial y_{0}}\left(  s\right)
\right\|  _{L^{\infty}\left(  y\right)  }ds\leq C\left\|  \Phi\right\|
_{Y_{0,\alpha}},\left\|  \frac{\partial\omega}{\partial y_{0}}\left(
\tau\right)  \right\|  _{L^{\infty}\left(  y\right)  }\leq Ce^{-\tau}\left\|
\Phi\right\|  _{Y_{0,\alpha}}, \label{omegay0}%
\end{equation}%
\begin{equation}
\int_{0}^{\tau}\left[  \frac{\partial\omega}{\partial y_{0}}\left(  s\right)
\right]  _{0,\alpha,\left(  y\right)  }ds\leq C\left\|  \Phi\right\|
_{Y_{\ell,\alpha}},\ \left[  \frac{\partial\omega}{\partial y_{0}}\left(
\tau\right)  \right]  _{0,\alpha,\left(  y\right)  }\leq Ce^{-\tau}\left\|
\Phi\right\|  _{Y_{\ell,\alpha}}. \label{omegay0a}%
\end{equation}

\end{lemma}

\begin{proof}
The method of proof is similar to the one used in the proof of Lemma
\ref{phidecayk0}.\textbf{\ }We take $\frac{\partial}{\partial y_{0}}$ of
(\ref{omega}) to get
\begin{align*}
\frac{\partial\omega}{\partial y_{0}}\left(  s\right)   &  =-\int_{s}^{\tau
}e^{-\left(  N-2\right)  \xi}\nabla_{y}^{2}\Phi\left(  \xi\right)
\{\frac{e^{-\xi}-e^{-\tau}}{1-e^{-\tau}}+\int_{0}^{\xi}e^{-\left(  \xi
-\eta\right)  }\frac{\partial\omega}{\partial y_{0}}\left(  \eta\right)
d\eta\}d\xi\\
&  +\frac{e^{-\tau}}{1-e^{-\tau}}\int_{0}^{\tau}e^{-\left(  N-3\right)  \xi
}\left(  1-e^{-\xi}\right)  \nabla_{y}^{2}\Phi\left(  Y\left(  \xi\right)
,\xi\right)  \{\frac{e^{-\xi}-e^{-\tau}}{1-e^{-\tau}}+\int_{0}^{\xi
}e^{-\left(  \xi-\eta\right)  }\frac{\partial\omega}{\partial y_{0}}\left(
\eta\right)  d\eta\}d\xi.
\end{align*}
Taking the $L^{\infty}\left(  y\right)  $ norm yields
\begin{align*}
\left\|  \frac{\partial\omega}{\partial y_{0}}\left(  s\right)  \right\|
_{L^{\infty}\left(  y\right)  }  &  \leq C\left\|  \Phi\right\|
_{Y_{0,\alpha}}\int_{s}^{\tau}e^{-\left(  N-2\right)  \xi}\{e^{-\xi}+e^{-\tau
}+\int_{0}^{\xi}e^{-\left(  \xi-\eta\right)  }\left\|  \frac{\partial\omega
}{\partial y_{0}}\left(  \eta\right)  \right\|  _{L^{\infty}\left(  y\right)
}d\eta\}d\xi\\
&  +C\left\|  \Phi\right\|  _{Y_{0,\alpha}}e^{-\tau}\int_{0}^{\tau}e^{-\left(
N-3\right)  \xi}\{e^{-\xi}+e^{-\tau}+\int_{0}^{\xi}e^{-\left(  \xi
-\eta\right)  }\left\|  \frac{\partial\omega}{\partial y_{0}}\left(
\eta\right)  \right\|  _{L^{\infty}\left(  y\right)  }d\eta\}d\xi.
\end{align*}
Integrating the above inequality from $s=0$ to $s=\tau$ and using
$e^{-(\xi-\eta)}\leq1,$ $e^{-\tau}\leq e^{-\xi}$ for $\eta\leq\xi$, $\xi
\leq\tau$ and $N\geq3$ yield
\begin{align*}
\int_{0}^{\tau}\left\|  \frac{\partial\omega}{\partial y_{0}}\left(  s\right)
\right\|  _{L^{\infty}\left(  y\right)  }ds  &  \leq C\left\|  \Phi\right\|
_{Y_{0,\alpha}}\int_{0}^{\tau}e^{-2s}ds+\left\|  \Phi\right\|  _{Y_{0,\alpha}%
}\left(  \int_{0}^{\tau}e^{-s}ds\right)  \left(  \int_{0}^{\tau}\left\|
\frac{\partial\omega}{\partial y_{0}}\left(  \eta\right)  \right\|
_{L^{\infty}\left(  y\right)  }d\eta\right) \\
&  +C\left\|  \Phi\right\|  _{Y_{0,\alpha}}e^{-\tau}\tau+C\left\|
\Phi\right\|  _{Y_{0,\alpha}}e^{-\tau}\tau^{2}\left(  \int_{0}^{\tau}\left\|
\frac{\partial\omega}{\partial y_{0}}\left(  \eta\right)  \right\|
_{L^{\infty}\left(  y\right)  }d\eta\right)  .
\end{align*}
Thus we have
\begin{equation}
\int_{0}^{\tau}\left\|  \frac{\partial\omega}{\partial y_{0}}\left(  s\right)
\right\|  _{L^{\infty}\left(  y\right)  }ds\leq C\left\|  \Phi\right\|
_{Y_{0,\alpha}}, \label{int-omega}%
\end{equation}
provided $\varepsilon_{0}$ is small enough. We now specialize to $s=\tau$ and
use (\ref{int-omega}) as well as the fact that $N\geq3$ to get
\begin{align*}
\left\|  \frac{\partial\omega}{\partial y_{0}}\left(  \tau\right)  \right\|
_{L^{\infty}\left(  y\right)  }  &  \leq C\left\|  \Phi\right\|
_{Y_{0,\alpha}}e^{-\tau}\int_{0}^{\tau}\{e^{-\xi}+\int_{0}^{\xi}e^{-\left(
\xi-\eta\right)  }\left\|  \frac{\partial\omega}{\partial y_{0}}\left(
\eta\right)  \right\|  _{L^{\infty}\left(  y\right)  }d\eta\}d\xi\\
&  \leq C\left\|  \Phi\right\|  _{Y_{0,\alpha}}e^{-\tau}+C\left\|
\Phi\right\|  _{Y_{0,\alpha}}e^{-\tau}\int_{0}^{\tau}\left\|  \frac
{\partial\omega}{\partial y_{0}}\left(  \eta\right)  \right\|  _{L^{\infty
}\left(  y\right)  }\left(  \int_{\eta}^{\tau}e^{-\left(  \xi-\eta\right)
}d\xi\right)  d\eta\\
&  \leq C\left\|  \Phi\right\|  _{Y_{0,\alpha}}e^{-\tau}%
\end{align*}
where we changed the order of integration. We thus obtain (\ref{omegay0}).
Using (\ref{omegay0}) we deduce (\ref{omegay0a}) in a similar way. Thus we
complete the proof.
\end{proof}

We also obtain the following estimates for the derivative of $\omega$ with
respect to $y$.

\begin{lemma}
\label{omegay}There exists $\varepsilon_{0}$ small such that for any $\tau
\geq1$ and any function $\Phi$ satisfying
\[
\left\|  \Phi\right\|  _{Y_{0,\alpha}}\leq\varepsilon_{0},
\]
we have
\[
\sup_{0\leq s\leq\tau}\left\|  \frac{\partial\omega}{\partial y}\left(
s\right)  \right\|  _{L^{\infty}\left(  y\right)  }\leq C\left\|
\Phi\right\|  _{Y_{0,\alpha}},\ \sup_{0\leq s\leq\tau}\left[  \frac
{\partial\omega}{\partial y}\left(  s\right)  \right]  _{0,\alpha,\left(
y\right)  }\leq C\left\|  \Phi\right\|  _{Y_{0,\alpha}}.
\]

\end{lemma}

\begin{proof}
We take $\frac{\partial}{\partial y}$ of (\ref{omega}) to get
\begin{align*}
\frac{\partial\omega}{\partial y}\left(  s\right)   &  =-\int_{s}^{\tau
}e^{-\left(  N-2\right)  \xi}\nabla_{y}^{2}\Phi\left(  Y\left(  \xi\right)
,\xi\right)  \{\frac{1-e^{-\xi}}{1-e^{-\tau}}+\int_{0}^{\xi}e^{-\left(
\xi-\eta\right)  }\frac{\partial\omega}{\partial y}\left(  \eta\right)
d\eta\}d\xi\\
&  +\frac{e^{-\tau}}{1-e^{-\tau}}\int_{0}^{\tau}e^{-\left(  N-3\right)  \xi
}\left(  1-e^{-\xi}\right)  \nabla_{y}^{2}\Phi\left(  Y\left(  \xi\right)
,\xi\right)  \{\frac{1-e^{-\xi}}{1-e^{-\tau}}+\int_{0}^{\xi}e^{-\left(
\xi-\eta\right)  }\frac{\partial\omega}{\partial y}\left(  \eta\right)
d\eta\}d\xi.
\end{align*}
Since $N\geq3$, we have
\begin{align*}
\sup_{0\leq s\leq\tau}\left\|  \frac{\partial\omega}{\partial y}\left(
s\right)  \right\|  _{L^{\infty}\left(  y\right)  }  &  \leq C\left\|
\Phi\right\|  _{Y_{0,\alpha}}\int_{s}^{\tau}e^{-\left(  N-2\right)  \xi
}\{1+\left(  \sup_{0\leq s\leq\tau}\left\|  \frac{\partial\omega}{\partial
y}\left(  s\right)  \right\|  _{L^{\infty}\left(  y\right)  }\right)  \int
_{0}^{\xi}e^{-\left(  \xi-\eta\right)  }d\eta\}d\xi\\
&  +C\left\|  \Phi\right\|  _{Y_{0,\alpha}}e^{-\tau}\int_{0}^{\tau}\{1+\left(
\sup_{0\leq s\leq\tau}\left\|  \frac{\partial\omega}{\partial y}\left(
s\right)  \right\|  _{L^{\infty}\left(  y\right)  }\right)  \int_{0}^{\xi
}e^{-\left(  \xi-\eta\right)  }d\eta\}d\xi\\
&  \leq C\left\|  \Phi\right\|  _{Y_{0,\alpha}}+C\left\|  \Phi\right\|
_{Y_{0,\alpha}}\left(  \sup_{0\leq s\leq\tau}\left\|  \frac{\partial\omega
}{\partial y}\left(  s\right)  \right\|  _{L^{\infty}\left(  y\right)
}\right)
\end{align*}
Thus we have
\begin{equation}
\sup_{0\leq s\leq\tau}\left\|  \frac{\partial\omega}{\partial y}\left(
s\right)  \right\|  _{L^{\infty}\left(  y\right)  }\leq C\left\|
\Phi\right\|  _{Y_{0,\alpha}}, \label{omega-inf}%
\end{equation}
provided $\varepsilon_{0}$ is small enough. In a similar manner, using
(\ref{omega-inf}) we obtain
\begin{align*}
\sup_{0\leq s\leq\tau}\left[  \frac{\partial\omega}{\partial y}\left(
s\right)  \right]  _{0,\alpha,\left(  y\right)  }  &  \leq C\left\|
\Phi\right\|  _{Y_{0,\alpha}}\int_{s}^{\tau}e^{-\left(  N-2\right)  \xi
}\{1+\left(  \sup_{0\leq s\leq\tau}\left[  \frac{\partial\omega}{\partial
y}\left(  s\right)  \right]  _{0,\alpha,\left(  y\right)  }\right)  \int
_{0}^{\xi}e^{-\left(  \xi-\eta\right)  }d\eta\}d\xi\\
&  +C\left\|  \Phi\right\|  _{Y_{0,\alpha}}e^{-\tau}\int_{0}^{\tau}\{1+\left(
\sup_{0\leq s\leq\tau}\left[  \frac{\partial\omega}{\partial y}\left(
s\right)  \right]  _{0,\alpha,\left(  y\right)  }\right)  \int_{0}^{\xi
}e^{-\left(  \xi-\eta\right)  }d\eta\}d\xi\\
&  \leq C\left\|  \Phi\right\|  _{Y_{0,\alpha}}+C\left\|  \Phi\right\|
_{Y_{0,\alpha}}\left(  \sup_{0\leq s\leq\tau}\left[  \frac{\partial\omega
}{\partial y}\left(  s\right)  \right]  _{0,\alpha,\left(  y\right)  }\right)
.
\end{align*}
This yields the H\"{o}lder estimate of $\frac{\partial\omega}{\partial y}$ and
completes the proof.
\end{proof}

We present the following estimates for higher-order derivatives similar to
Theorems in the previous section.

\begin{lemma}
\label{omegal}Let $\ell\geq1$ be an integer. There exists $\varepsilon_{0}$
small such that for any $\tau\geq1$ and any function $\Phi$ satisfying
\[
\left\|  \Phi\right\|  _{Y_{\ell,\alpha}}\leq\varepsilon_{0},
\]
we have the following
\[
\sup_{0\leq s\leq\tau}\left\|  \frac{\partial^{\ell}\omega}{\partial y^{\ell}%
}\left(  s\right)  \right\|  _{L^{\infty}\left(  y\right)  }\leq C\left\|
\Phi\right\|  _{Y_{\ell,\alpha}},\ \sup_{0\leq s\leq\tau}\left[
\frac{\partial^{\ell}\omega}{\partial y^{\ell}}\left(  s\right)  \right]
_{0,\alpha,\left(  y\right)  }\leq C\left\|  \Phi\right\|  _{Y_{\ell,\alpha}%
},
\]%
\[
\int_{0}^{\tau}\left\|  \frac{\partial^{\ell+1}\omega}{\partial y^{\ell
}\partial y_{0}}\left(  s\right)  \right\|  _{L^{\infty}\left(  y\right)
}ds\leq C\left\|  \Phi\right\|  _{Y_{\ell,\alpha}},\ \int_{0}^{\tau}\left[
\frac{\partial^{\ell+1}\omega}{\partial y^{\ell}\partial y_{0}}\left(
s\right)  \right]  _{0,\alpha,\left(  y\right)  }ds\leq C\left\|
\Phi\right\|  _{Y_{\ell,\alpha}},
\]%
\[
\ \left\|  \frac{\partial^{\ell+1}\omega}{\partial y^{\ell}\partial y_{0}%
}\left(  \tau\right)  \right\|  _{L^{\infty}\left(  y\right)  }\leq Ce^{-\tau
}\left\|  \Phi\right\|  _{Y_{\ell,\alpha}},\left[  \frac{\partial^{\ell
+1}\omega}{\partial y^{\ell}\partial y_{0}}\left(  \tau\right)  \right]
_{0,\alpha,\left(  y\right)  }\leq Ce^{-\tau}\left\|  \Phi\right\|  _{Y_{\ell
}}.
\]

\end{lemma}

\bigskip{}

As a consequence of Theorem \ref{thmSS}, we have, for any $k\geq0$ integer and
$0<\alpha<1$,
\begin{equation}
\left\|  \rho\right\|  _{C^{k,\alpha}\left(  \mathbb{R}^{N}\right)  }\leq
C\varepsilon_{0}\;,\ \left\|  \Phi\right\|  _{C^{k+2,\alpha}\left(
\mathbb{R}^{N}\right)  }\leq C\varepsilon_{0}. \label{rhoPhi}%
\end{equation}
We now study the limit behaviour of the self-similar system (\ref{SS1}%
)-(\ref{SS3}). Indeed, the limit behaviour is asymptotically equivalent to the
free streaming case.

\subsection{Proof of Theorem \ref{asymp}.}

\begin{proof}
We begin with
\begin{align}
\omega\left(  s;y,y_{0},\tau\right)   &  =-\int_{s}^{\tau}e^{-\left(
N-2\right)  \xi}\nabla_{y}\Phi\left(  Y\left(  \xi;y,y_{0},\tau\right)
,\xi\right)  d\xi\label{omega rep}\\
&  +\frac{e^{-\tau}}{1-e^{-\tau}}\int_{0}^{\tau}e^{-\left(  N-3\right)  \xi
}\left(  1-e^{-\xi}\right)  \nabla_{y}\Phi\left(  Y\left(  \xi;y,y_{0}%
,\tau\right)  ,\xi\right)  d\xi.\nonumber
\end{align}
By using (\ref{rhoPhi}) and by the dominated convergence theorem, as
$\tau\rightarrow\infty,$ in $C^{k+1},$%
\[
\omega\left(  s;y,y_{0},\tau\right)  \rightarrow\omega\left(  s;y,y_{0}%
,\infty\right)  \equiv\omega_{\infty}\left(  s,y,y_{0}\right)  .
\]
In particular, we have, as $\tau\rightarrow\infty,$ in $C^{k+1},$%
\[
v_{0}\left(  y,y_{0},\tau\right)  =v_{0}\left(  y,v\left(  y,y_{0}%
,\tau\right)  \right)  =\tilde{V}+\omega\left(  0;y,y_{0},\tau\right)
\rightarrow y+\omega_{\infty}\left(  0,y,y_{0}\right)  .
\]
Thus, we have, as $\tau\rightarrow\infty,$ in $C_{loc}^{k},$%
\[
e^{-N\tau}g\left(  y,y_{0},\tau\right)  =g_{0}\left(  y_{0},v_{0}\left(
y,y_{0},\tau\right)  \right)  \rightarrow g_{0}\left(  y_{0},y+\omega_{\infty
}\left(  0,y,y_{0}\right)  \right)  \equiv g_{\infty}\left(  y,y_{0}\right)  .
\]
Next, since
\[
e^{\tau}\frac{\partial v}{\partial y_{0}}=-\left(  \frac{1}{1-e^{-\tau}%
}\right)  I_{N}+e^{\tau}\frac{\partial\omega}{\partial y_{0}}\left(  t\right)
,
\]
using Lemma \ref{omega0}-Lemma \ref{omegal} yields, as $\tau\rightarrow
\infty,$ in $C^{k},$
\[
e^{N\tau}\left\vert \det\frac{\partial v}{\partial y_{0}}\right\vert
\rightarrow J_{\infty}\left(  y,y_{0}\right)  \simeq1+\mathcal{O}\left(
\varepsilon_{0}\right)  .
\]
We then apply the dominated convergence theorem to get, as $\tau
\rightarrow\infty,$ in $C_{loc}^{k},$%
\begin{align*}
\bar{\rho}\left(  y,\tau\right)   &  =\int g_{0}\left(  y_{0},y+\omega\left(
0;y,y_{0},\tau\right)  \right)  e^{N\tau}\det\frac{\partial v}{\partial y_{0}%
}dy_{0}\\
&  \rightarrow\int g_{0}\left(  y_{0},y+\omega_{\infty}\left(  0,y,y_{0}%
\right)  \right)  J_{\infty}\left(  y,y_{0}\right)  dy_{0}\\
&  \equiv\bar{\rho}_{\infty}\left(  y\right)  .
\end{align*}
Using the elliptic regularity theory from the equation
\[
\Delta_{y}\Phi=\bar{\rho},
\]
there exists $\Phi_{\infty}\left(  y\right)  \in C_{loc}^{k+1,\beta}$ , for
any $0<\beta<1,$ such that
\[
\Delta_{y}\Phi_{\infty}=\bar{\rho}_{\infty}.
\]
Taking the limit in (\ref{g}) as $\tau\rightarrow\infty$ yields 
(\ref{asymp g}). Finally, notice that, given a test 
function\textbf{\ }$\psi\left(  y,v\right)
$%

\begin{align*}
\int_{\mathbb{R}^{3}\times\mathbb{R}^{3}}g\left(  y,v,\tau\right)  \psi\left(
y,v\right)  dydv &  =\int e^{N\tau}g_{0}\left(  y_{0},v_{0}\right)
\psi\left(  y,v\right)  dydv=\\
&  \int e^{N\tau}g_{0}\left(  y_{0},w_{0}\left(  y,y_{0},\tau\right)  \right)
\psi\left(  y,w\left(  y,y_{0},\tau\right)  \right)  \det\frac{\partial
v}{\partial y_{0}}dydy_{0}%
\end{align*}
and taking the limit $\tau\rightarrow\infty$ we obtain
\[
\int_{\mathbb{R}^{3}\times\mathbb{R}^{3}}g\left(  y,v,\tau\right)  \psi\left(
y,v\right)  dydv\rightarrow\int g_{0}\left(  y_{0},y+\omega_{\infty}\left(
0,y,y_{0}\right)  \right)  \psi\left(  y,y\right)  J_{\infty}\left(
y,y_{0}\right)  dy_{0}dy
\]
which can be written in the sense of distributions as
\[
g\left(  y,v,\tau\right)  \rightarrow\int g_{\infty}\left(  y_{0},y\right)
\delta\left(  v-y\right)  J_{\infty}\left(  y,y_{0}\right)  dy_{0}%
\;\;\text{as\ \ }\tau\rightarrow\infty
\]

This yields (\ref{S4E1}), whence the proof is complete.

Notice that in the limit case $\varepsilon_{0}\rightarrow0$ (\ref{S4E1})
reduces to
\begin{align*}
g\left(  y,v,\tau\right)   &  \rightarrow\left[  \int g_{\infty}\left(
y_{0},y\right)  dy_{0}\right]  \delta\left(  y-v\right)  \\
&  =\left[  \int g_{0}\left(  y_{0},y\right)  dy_{0}\right]  \delta\left(
y-v\right)  \;\;\text{as\ \ }\tau\rightarrow\infty
\end{align*}

\end{proof}

\bigskip

\textbf{Acknowledgements: }Hyung Ju\ Hwang and Juan Vel\'{a}zquez acknowledge
the hospitality and support of the Max Planck Institute for Gravitational
Physics where this work was begun. Hyung Ju Hwang is supported by 
SFI Grant \ 06/RFP/MAT029. Juan Vel\'{a}zquez acknowledges also the
support of the Humboldt Foundation, the Max Planck Institute for Mathematics
in the Natural Sciences and DGES Grant MTM2004-05634.

\bigskip

\bigskip

\bigskip

\bigskip

\bigskip

\bigskip

\end{document}